\documentclass[12pt]{article}
\usepackage{amsfonts}
\usepackage{amsmath}
\usepackage{amssymb}
\usepackage[english]{babel}
\usepackage[T1]{fontenc}
\usepackage{graphicx}
\usepackage{subfigure}
\newtheorem{prop}{Proposition}[section]
\newtheorem{rem}[prop]{Remark}

\newtheorem{corr}[prop]{Corollary}

\usepackage{epsfig}
\usepackage{latexsym}
\usepackage{amstext}
\usepackage{epsfig, graphicx}
\usepackage{epsfig}
\usepackage{latexsym}
\usepackage{amstext}
\usepackage{amssymb}
\usepackage{graphics}
\usepackage{enumerate}

\newcommand{\beq}{\begin{eqnarray}}
\newcommand{\beqq}{\begin{eqnarray*}}
\newcommand{\eeq}{\end{eqnarray}}
\newcommand{\eeqq}{\end{eqnarray*}}

\title{On hitting times of the winding processes of planar Brownian motion and of Ornstein-Uhlenbeck processes, via Bougerol's identity.}
\author{ S. Vakeroudis \thanks{Laboratoire de Probabilit\'{e}s et Mod\`{e}les
Al\'{e}atoires (LPMA) CNRS : UMR7599,  Universit\'{e} Pierre et Marie
Curie - Paris VI, 4 Place Jussieu, 75252 Paris Cedex 05, France.
E-mail: stavros.vakeroudis@etu.upmc.fr} \thanks{D\'{e}partement de Math\'{e}matiques et de Biologie CNRS : UMR8544,  Ecole Normale Sup\'{e}rieure, 46 rue d'Ulm, 75005 Paris, France.} }
\date{\today}
\begin{document}
%%%%%%%%%%%%%%%%%%%%%%%%%%%%%%%%%%%%%%%%%%%%%%%%%%%%%%%%%%%%%%%%

\maketitle

\begin{abstract}
Some identities in law in terms of planar complex valued Ornstein-Uhlenbeck processes $(Z_{t}=X_{t}+iY_{t},t\geq0)$ including planar Brownian motion are established and shown to be equivalent to the well known Bougerol identity for linear Brownian motion $(\beta_{t},t\geq0)$: for any fixed $u>0$:
    \beqq \sinh(\beta_{u}) \stackrel{(law)}{=}
    \hat{\beta}_{(\int^{u}_{0}ds\exp(2\beta_{s}))},
\eeqq
with $(\hat{\beta}_{t},t\geq0)$ a Brownian motion, independent of $\beta$. \\
These identities in law for 2-dimensional processes allow to study the distributions of hitting times
$T^{\theta}_{c}\equiv\inf\{ t:\theta_{t} =c \}, \ (c>0)$, \,
$T^{\theta}_{-d,c}\equiv\inf\{ t:\theta_{t}\notin(-d,c) \}, \
(c,d>0)$ \ and more specifically of $T^{\theta}_{-c,c}\equiv\inf\{
t:\theta_{t}\notin(-c,c) \}, \ (c>0)$ \ of the continuous winding processes \
$\theta_{t}=\mathrm{Im}(\int^{t}_{0}\frac{dZ_{s}}{Z_{s}}), t\geq0$ of complex valued Ornstein-Uhlenbeck processes.
\end{abstract}

%%%%%%%%%%%%%%%%%%%%%%%%%%%%%%%%%%%%%%%%%%%%%%%%%%%%%%%%%%%%%%%%
$\vspace{5pt}$

%%%%%%%%%%%%%%%%%%%%%%%%%%%%%%%%%%%%%%%%%%%%%%%%%%%%%%%%%%%%%%%%
\section{Introduction}
%%%%%%%%%%%%%%%%%%%%%%%%%%%%%%%%%%%%%%%%%%%%%%%%%%%%%%%%%%%%%%%%

\renewcommand{\thefootnote}{\fnsymbol{footnote}}
The conformal invariance of planar Brownian motion has deep consequences as to the structure of its trajectories (see, e.g., Le Gall \cite{LeG90}). In particular, a number of articles have been devoted to the study of its continuous winding process $(\theta _{t},t\geq0)$: Spitzer \cite{Spi58}, Williams \cite{Wil74},  Durrett \cite{Dur82}, Messulam-Yor \cite{MeY82}, Pitman-Yor \cite{PiY86}, Le Gall-Yor \cite{LeGY87}, Bertoin-Werner \cite{BeW94}, Yor \cite{Yor97}, Pap-Yor \cite{PaY00}, Bentkus-Pap-Yor \cite{BPY03}. In this paper, we take up again the study of the first hitting times:
    \beqq T^{\theta}_{-d,c}\equiv\inf\{ t:\theta_{t}\notin(-d,c) \}, \ (c,d>0),
\eeqq
this time in relation with Bougerol's well-known identity (see Bougerol \cite{Bou83}, Alili-Dufresne-Yor \cite{ADY97} and Yor \cite{Yor01}): for fixed $u>0$:
    \beqq \sinh(\beta_{u}) \stackrel{(law)}{=} \hat{\beta}_{(\int^{u}_{0}ds\exp(2\beta_{s}))} \ ,
\eeqq
where $(\hat{\beta_{t}},t\geq0)$ is a Brownian motion\footnote[3]{When we simply write: Brownian motion, we always mean real-valued Brownian motion, starting from 0. For 2-dimensional Brownian motion, we indicate planar or complex BM.}, independent of $\beta$. \\
In particular, it turns out that: for fixed $c>0$:
\beqq
\theta_{T^{\hat{\beta}}_{c}} \stackrel{(law)}{=} C_{a(c)}, \ \ \ \ \ \ \ (\star)
\eeqq
where $\hat{\beta}$ is a BM$^{\ddag}$ independent of $(\theta _{u},u\geq0)$, $T^{\hat{\beta}}_{c}=\inf\{t:\hat{\beta_{t}}=c \}$,
$(C_{t},t\geq0)$ is a standard Cauchy process and $a(c)= \arg \sinh (c)\equiv \log \left(c+\sqrt{1+c^{2}}\right), \; c\in \mathbb{R}$. \\
The identity $(\star)$ yields yet another proof of the celebrated Spitzer theorem:
\beqq
\frac{2}{\log t} \; \theta_{t} \overset{{(law)}}{\underset{t\rightarrow\infty}\longrightarrow} C_{1} \ ,
\eeqq
with the help of Williams' "pinching method" (see Williams \cite{Wil74} and Messulam-Yor \cite{MeY82}).

Moreover, we study the distributions of $ T^{\theta}_{-\infty,c}$ and $T^{\theta}_{-c,c}$. In particular, we give explicit formulae for the density function of $T^{\theta}_{-c,c}$ and for the first moment of $\ln \left(T^{\theta}_{-c,c}\right)$.

The last section of the paper is devoted to developing similar results when planar Brownian motion is replaced by a complex valued Ornstein-Uhlenbeck process. We note that Bertoin-Werner \cite{BeW94} already made discussions of windings for planar Brownian motion using arguments related to Ornstein-Uhlenbeck processes.

Firstly, we obtain some analogue of $(\star)$ when $T^{\hat{\beta}}_{c}$ is replaced by $T^{(\lambda)}_{c}=T^{\theta^{Z}}_{-c,c}=\inf\{t:|\theta^{Z}_{t}|=c \}$, the corresponding time for an Ornstein-Uhlenbeck process with parameter $\lambda$. Secondly, we identify the distribution of $T^{(\lambda)}_{c}$. More specifically, we derive the asymptotics of $E\left[T^{(\lambda)}_{c}\right]$ for $\lambda$ large and for $\lambda$ small.

%%%%%%%%%%%%%%%%%%%%%%%%%%%%%%%%%%%%%%%%%%%%%%%%%%%%%%%%%%%%%%%%
\section{The Brownian motion case}
%%%%%%%%%%%%%%%%%%%%%%%%%%%%%%%%%%%%%%%%%%%%%%%%%%%%%%%%%%%%%%%%

%%%%%%%%%%%%%%%%%%%%%%%%%%%%%%%%%%%%%%%%%%%%%%%%%%%%%%%%%%%%%%%%
\subsection{A reminder on planar Brownian motion}
%%%%%%%%%%%%%%%%%%%%%%%%%%%%%%%%%%%%%%%%%%%%%%%%%%%%%%%%%%%%%%%%

Let $(Z_{t}=X_{t}+iY_{t},t\geq0)$ denote a standard planar Brownian motion, starting from $x_{0}+i0,x_{0}>0$, where $(X_{t},t\geq0)$ and $(Y_{t},t\geq0)$ are two independent linear Brownian motions, starting respectively from $x_{0}$ and $0$.

As is well known (see e.g. It\^{o}-McKean \cite{ItMK65}), since $x_{0}\neq0$,
$(Z_{t},t\geq0)$ does not visit a.s. the point $0$ but keeps
winding around $0$ infinitely often. In particular,
the continuous winding process $\theta_{t}=\mathrm{Im}(\int^{t}_{0}\frac{dZ_{s}}{Z_{s}}),t\geq0$ is well
defined.

Furthermore, there is the skew product representation:
\beq\label{skew-product}
\log\left|Z_{t}\right|+i\theta_{t}\equiv\int^{t}_{0}\frac{dZ_{s}}{Z_{s}}=\left(
\beta_{u}+i\gamma_{u}\right)
\Bigm|_{u=H_{t}=\int^{t}_{0}\frac{ds}{\left|Z_{s}\right|^{2}}},
\eeq
where $(\beta_{u}+i\gamma_{u},u\geq0)$ is another planar Brownian motion starting from $\log x_{0}+i0$. For a study of the Bessel clock $H$, see Yor \cite{Yor80}.

Rewriting (\ref{skew-product}) as:
\beq\label{skew-product2}
\log\left|Z_{t}\right|=\beta_{H_{t}}; \ \ \theta_{t}=\gamma_{H_{t}},
\eeq
we easily obtain that the total $\sigma$-fields $\sigma \{\left|Z_{t}\right|,t\geq0\}$ and $\sigma \{\beta_{u},u\geq0\}$ are identical, whereas $(\gamma_{u},u\geq0)$ is independent from $(\left|Z_{t}\right|,t\geq0)$. \\

%%%%%%%%%%%%%%%%%%%%%%%%%%%%%%%%%%%%%%%%%%%%%%%%%%%%%%%%%%%%%%%%
\begin{figure}
    \centering
        \includegraphics[width=1.00\textwidth]{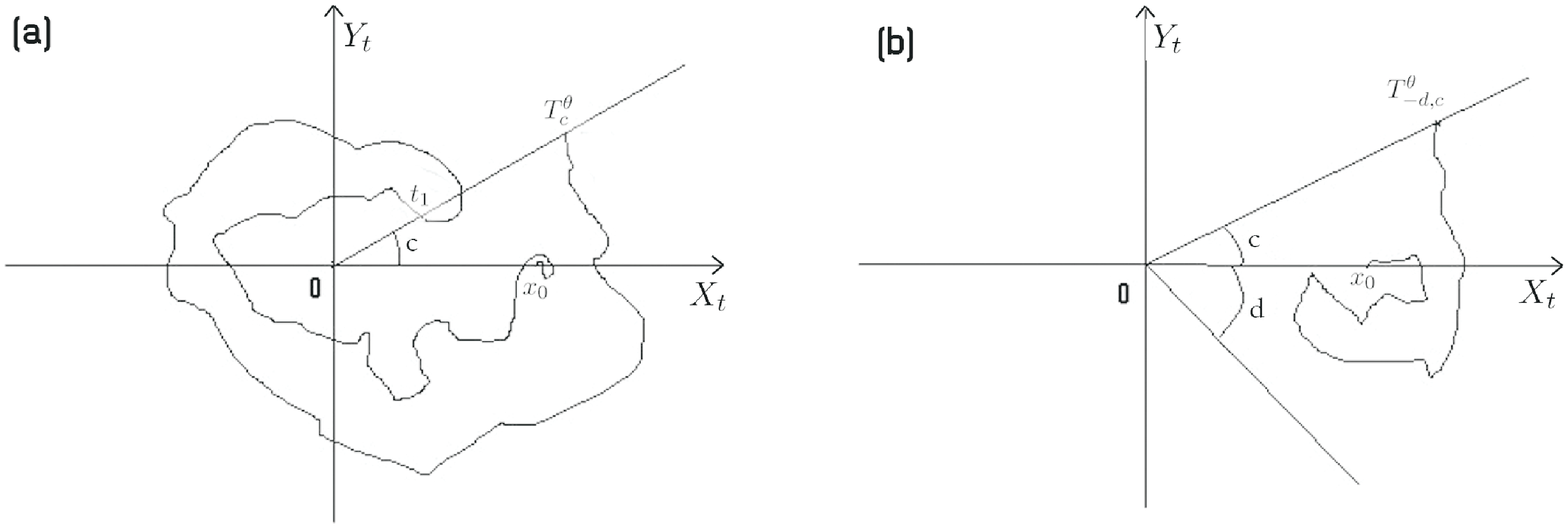}
        \caption{Exit times for a planar BM. This figure presents the exit times (a) $T^{\theta}_{c}$ ($t_{1}$ doesn't matter because the angle is negative)  and (b) $T^{\theta}_{-d,c}$ for a planar BM starting from $x_{0}+i0$. }\label{figexittimes}
\end{figure}
%%%%%%%%%%%%%%%%%%%%%%%%%%%%%%%%%%%%%%%%%%%%%%%%%%%%%%%%%%%%%%%%

A number of studies of the properties of the first hitting time  (see Figure \ref{figexittimes}(b))
\beqq
T^{\theta}_{-d,c}\equiv\inf\{ t:\theta_{t}\notin(-d,c) \},\ \ (c,d>0),
\eeqq
have been developed, going back to Spitzer \cite{Spi58}.\\
In particular, it is well known (Spitzer \cite{Spi58}, Burkholder \cite{Bur77},
Revuz-Yor \cite{ReY99} Ex. 2.21/page 196) that:
\beq\label{spi}
E\left[ (T^{\theta}_{-d,c})^{p} \right]<\infty \ \ \ \;
\mathrm{if}
\; \mathrm{and} \; \mathrm{only} \; \mathrm{if} \ \ \ \;
p<\frac{\pi}{2(c+d)}.
\eeq
Moreover, Spitzer's asymptotic theorem (see e.g. Spitzer \cite{Spi58}) states that:
\beq\label{spi2}
  \frac{2 \theta_{t}}{\log t} \overset{{(law)}}{\underset{t\rightarrow\infty}\longrightarrow} C_{1} \stackrel{(law)}{=} \gamma_{T_{1}^{\beta}},
\eeq
where $C_{1}$ is a standard Cauchy variable.\\

%%%%%%%%%%%%%%%%%%%%%%%%%%%%%%%%%%%%%%%%%%%%%%%%%%%%%%%%%%%%%%%%
\subsection{On the Laplace transform of the distribution of the hitting time $T^{\theta}_{c}\equiv T^{\theta}_{-\infty,c}$}
%%%%%%%%%%%%%%%%%%%%%%%%%%%%%%%%%%%%%%%%%%%%%%%%%%%%%%%%%%%%%%%%
Now, we use the representation (\ref{skew-product2}) to access the
distribution of $T^{\theta}_{c}$ (see Figure
\ref{figexittimes}(a)). We define $T^{\gamma}_{c}\equiv\inf\{
t:\gamma_{t}\notin(-\infty,c) \}$
the hitting time associated to the Brownian motion $(\gamma_{t}, t\geq 0)$.
Note that, from (\ref{skew-product2}): \\
$H_{T^{\theta}_{c}}=T^{\gamma}_{c}$, \
hence: $T^{\theta}_{c}=H^{-1}_{u \bigm|_{u=T^{\gamma}_{c}}}$, where
\beq\label{Hinverse}
H^{-1}_{u}\equiv \inf\{ t:H_{t}>u \} =
\int^{u}_{0}ds \exp(2\beta_{s}) := A_{u} .
\eeq
Thus, we have obtained:
\beq\label{TcA}
    T^{\theta}_{c}=A_{T^{\gamma}_{c}},
\eeq
where $(A_{u},u\geq0)$ and $T^{\gamma}_{c}$ are independent, since $\beta$ and $\gamma$ are independent. \\
We can write: $\beta_{s}=(\log x_{0})+ \beta_{s}^{(0)}$, with
$(\beta_{s}^{(0)},s\geq0)$ a standard one-dimensional Brownian
motion starting from $0$. Then, we deduce from (\ref{TcA}) that:
\beq\label{z0dep}
T^{\theta}_{c}=x_{0}^{2} \left(
\int^{T^{\gamma}_{c}}_{0}ds \exp(2\beta_{s}^{(0)}) \right).
\eeq
From now on, for simplicity, we shall take $x_{0}=1$, but this is really no restriction, as the dependency in $x_{0}$, which is exhibited in (\ref{z0dep}), is very simple.\\
We shall also make use of Bougerol's identity
\cite{Bou83,ADY97} and \cite{Yor01} (p. 200), which is very useful to study the
distribution of $A_{u}$ (e.g. \cite{MaY98,MaY05}). For any fixed $u>0$:
\beq\label{bougerolidentity}
    \sinh(\beta_{u}) \stackrel{(law)}{=}  \hat{\beta}_{A_{u}} = \hat{\beta}_{(\int^{u}_{0}ds\exp(2\beta_{s}))} ,
\eeq
where on the right hand side, $(\hat{\beta_{t}},t\geq0)$ is a
Brownian motion, independent of $A_{u}\equiv
\int^{u}_{0}ds\exp(2\beta_{s})$.
\\
Thus, from (\ref{bougerolidentity}) and (\ref{TcA}), and as is
well known \cite{ReY99}, the law of
$\beta_{T^{\gamma}_{c}}$ is the Cauchy law with parameter $c$,
i.e., with density:
\beqq
    h_{c}(y)=\frac{c}{\pi(c^{2}+y^{2})},
\eeqq
we deduce that:
%%%%%%%%%%%%%%%%%%%%%%%%%%%%%%%%%%%%%%%%%%%%%%%%%%%%%%%%%%%%%%%%
\begin{prop}\label{bougerol}
For fixed $c>0$, there is the following identity in law:
\beq\label{bougerolidentity2}
\sinh(C_{c}) \stackrel{(law)}{=}
\hat{\beta}_{(T^{\theta}_{c})},
\eeq
where, on the left hand side,
$(C_{c},c\geq0)$ denotes a standard Cauchy process and on the
right hand side, $(\hat{\beta}_{u},u\geq0)$ is a one-dimensional
BM, independent from $T^{\theta}_{c}$.
\end{prop}
%%%%%%%%%%%%%%%%%%%%%%%%%%%%%%%%%%%%%%%%%%%%%%%%%%%%%%%%%%%%%%%%
We may now identify the densities of the variables found on both sides of (\ref{bougerolidentity2}),i.e.: \\
on the left hand side: $\frac{1}{\sqrt{1+x^{2}}} h_{c}(\arg \sinh x) = \frac{1}{\sqrt{1+x^{2}}} h_{c}(a(x))$ ; \\
on the right hand side: $E \left[ \frac{1}{\sqrt{2\pi
T^{\theta}_{c}}} \exp \left( -\frac{x^{2}}{2T^{\theta}_{c}} \right) \right]$, \\
where $a(x)=\arg \sinh (x)$. \\
\vspace{0.3cm}
Thus, we have obtained the following:
%%%%%%%%%%%%%%%%%%%%%%%%%%%%%%%%%%%%%%%%%%%%%%%%%%%%%%%%%%%%%%%%
\begin{prop}\label{laplace}
The distribution of $T^{\theta}_{c}$ may be characterized by:
\beq\label{laplacetransform}
    E \left[ \frac{1}{\sqrt{2\pi T^{\theta}_{c}}} \exp \left( -\frac{x}{2T^{\theta}_{c}} \right) \right] = \frac{1}{\sqrt{1+x}} \: \frac{c}{\pi(c^{2}+\log^{2}(\sqrt{x}+\sqrt{1+x}))}, \ \ x \geq 0. \nonumber \\
\eeq
\end{prop}
%%%%%%%%%%%%%%%%%%%%%%%%%%%%%%%%%%%%%%%%%%%%%%%%%%%%%%%%%%%%%%%%
The proof of Proposition \ref{laplace} follows from: $a(y) = \arg
\sinh (y) \equiv \log(y+\sqrt{1+y^{2}})$ and by making the change
of variable $y^{2}=x$. Let us now define the probability:
\beqq
    Q_{c} =  \sqrt{ \frac{\pi c^{2}}{2T^{\theta}_{c}}} \cdot P \ .
\eeqq
The fact that $Q_{c}$ is a probability follows from
(\ref{laplacetransform}) by taking $x=0$. Thus we obtain that $c
\: E[\sqrt{\pi/2 T^{\theta}_{c}}]=1$, and we may write:
\beq\label{Qlaplacetransform}
    E_{Q_{c}} \left[ \exp \left( -\frac{x}{2T^{\theta}_{c}} \right) \right] = \frac{1}{\sqrt{1+x}}\frac{1}{1+\frac{1}{c^{2}}\log^{2}(\sqrt{x}+\sqrt{1+x})}, \ \ \forall x \geq 0,
\eeq
which yields the Laplace transform of $1/T^{\theta}_{c}$
under $Q_{c}$. \\
Let us now take a look at what happens if we make
$c\rightarrow\infty$. If we denote by $T_{1}^{\beta} \equiv \inf
\{ t: \beta_{t}=1 \}$ the first hitting time of level 1 for a
standard BM $\beta$ and by $N$ a standard Gaussian variable
$\mathcal{N}(0,1)$, from equation (\ref{Qlaplacetransform}), we
obtain:
\beq\label{Qlaplacetransform2}
    \underset{{c\rightarrow\infty}}{\lim }E_{Q_{c}} \left[ e^{-x/2T^{\theta}_{c}} \right] = E \left( e^{-x N^{2}/2} \right) = E \left( e^{-x/2T_{1}^{\beta}} \right),
\eeq
which means that :
$T^{\theta}_{c}\overset{{(law)}}{\underset{c\rightarrow\infty}\longrightarrow}T_{1}^{\beta}$.
(At this point, one may wonder whether there is some kind of
convergence in law involving $(\theta_{u},u\geq0)$, under $Q_{c}$,
as $c\rightarrow\infty$, but, we shall not touch this point).
\vspace{10pt}
\\
From Proposition \ref{laplace} we deduce the following:
%%%%%%%%%%%%%%%%%%%%%%%%%%%%%%%%%%%%%%%%%%%%%%%%%%%%%%%%%%%%%%%%
\begin{corr}\label{laplace2} Let $\varphi(x)$ denote the Laplace transform (\ref{Qlaplacetransform}), that is the Laplace transform of $1/2T^{\theta}_{c}$ under $Q_{c}$. Then, the Laplace transform of $1/2T^{\theta}_{c}$ under $P$ is:
\beq\label{laplacetransform2}
    E \left[ \exp \left( -\frac{x}{2T^{\theta}_{c}} \right) \right] = \int^{\infty}_{x} \frac{dw}{\sqrt{w-x}} \: \varphi(w).
\eeq
\end{corr}
%%%%%%%%%%%%%%%%%%%%%%%%%%%%%%%%%%%%%%%%%%%%%%%%%%%%%%%%%%%%%%%%
%%%%%%%%%%%%%%%%%%%%%%%%%%%%%%%%%%%%%%%%%%%%%%%%%%%%%%%%%%%%%%%%
\textbf{Proof of Corollary \ref{laplace2}} From Fubini's theorem,
we deduce from (\ref{Qlaplacetransform}) that:
\beqq
    E \left[ \exp \left( -\frac{x}{2T^{\theta}_{c}} \right) \right] &=& \int^{\infty}_{0} \frac{dy}{\sqrt{y}} \: E \left[ \frac{1}{\sqrt{2\pi T^{\theta}_{c}}} \exp \left( -\frac{x+y}{2T^{\theta}_{c}} \right) \right] \\
    &=& \int^{\infty}_{0} \frac{dy}{\sqrt{y}} \: \varphi(x+y) \\
    &\stackrel{y=xt}{=}& \sqrt{x} \int^{\infty}_{0} \frac{dt}{\sqrt{t}} \: \varphi(x(1+t)) \\
    &\stackrel{v=1+t}{=}& \sqrt{x} \int^{\infty}_{1} \frac{dv}{\sqrt{v-1}} \: \varphi(xv) \\
    &\stackrel{w=xv}{=}& \int^{\infty}_{x} \frac{dw}{\sqrt{w-x}} \: \varphi(w),
\eeqq
which is formula (\ref{laplacetransform2}).
\begin{flushright}
  $\Box$
\end{flushright}
%%%%%%%%%%%%%%%%%%%%%%%%%%%%%%%%%%%%%%%%%%%%%%%%%%%%%%%%%%%%%%%%

%%%%%%%%%%%%%%%%%%%%%%%%%%%%%%%%%%%%%%%%%%%%%%%%%%%%%%%%%%%%%%%%
\subsection{Some related identities in law}
%%%%%%%%%%%%%%%%%%%%%%%%%%%%%%%%%%%%%%%%%%%%%%%%%%%%%%%%%%%%%%%%
This subsection is strongly related to \cite{DuY10}. \\
A slightly different look at the combination of Bougerol's
identity (\ref{bougerolidentity}) and the skew-product
representation (\ref{skew-product}) lead to the following
striking identities in law:
%%%%%%%%%%%%%%%%%%%%%%%%%%%%%%%%%%%%%%%%%%%%%%%%%%%%%%%%%%%%%%%%
\begin{prop}\label{law} Let $(\delta_{u},u\geq0)$ be a 1-dimensional Brownian
motion independent of the planar Brownian motion $(Z_{u},u\geq
0)$, starting from $1+i0$. Then, for any $b\geq0$, the following identities in law hold:
\beqq
    \mathrm{(i)} \: H_{T^{\delta}_{b}} \stackrel{(law)}{=} T^{\beta}_{a(b)} \ \ \ \
    \mathrm{(ii)} \: \theta_{T^{\delta}_{b}} \stackrel{(law)}{=} C_{a(b)} \ \ \ \
    \mathrm{(iii)} \: \bar{\theta}_{T^{\delta}_{b}} \stackrel{(law)}{=}
    |C_{a(b)}|,
\eeqq
where $C_{A}$ is a Cauchy variable with parameter $A$ and
$\bar{\theta}_{u}=\sup_{s\leq u}\theta_{s}$.
\end{prop}
%%%%%%%%%%%%%%%%%%%%%%%%%%%%%%%%%%%%%%%%%%%%%%%%%%%%%%%%%%%%%%%%
%%%%%%%%%%%%%%%%%%%%%%%%%%%%%%%%%%%%%%%%%%%%%%%%%%%%%%%%%%%%%%%%
\textbf{Proof of Proposition \ref{law}} From the symmetry
principle (see \cite{And87} for the original Note and \cite{Gal08} for a detailed discussion), Bougerol's
identity may be equivalently stated as:
\beq\label{bougerol2}
    \sinh(\bar{\beta}_{u}) \stackrel{(law)}{=} \bar{\delta}_{A_{u}(\beta)}.
\eeq
Consequently, the laws of the first hitting times of a fixed
level $b$ by the processes on each side of (\ref{bougerol2}) are
identical, that is:
\beqq
T^{\beta}_{a(b)} \stackrel{(law)}{=} H_{T^{\delta}_{b}},
\eeqq
which is (i). \\
(ii) follows from (i) since:
\beqq
\theta_{u} \stackrel{(law)}{=} \gamma_{H_{u}},
\eeqq
with $(\gamma_{s},s\geq0)$ a Brownian motion independent of
$(H_{u},u\geq0)$ and $(C_{u},u\geq0)$ may be represented as
$(\gamma_{T^{\beta}_{u}},u\geq0)$. \\
(iii) follows from (ii), again with the help of the symmetry
principle.
\begin{flushright}
  $\Box$
\end{flushright}
%%%%%%%%%%%%%%%%%%%%%%%%%%%%%%%%%%%%%%%%%%%%%%%%%%%%%%%%%%%%%%%%
%%%%%%%%%%%%%%%%%%%%%%%%%%%%%%%%%%%%%%%%%%%%%%%%%%%%%%%%%%%%%%%%
\begin{rem}
Proposition \ref{laplace} may also be
derived from (iii) in Proposition \ref{law}. Indeed, for $c>0$, starting from
the LHS of (iii), and letting $N\thicksim \mathcal{N}(0,1)$
independent from $T^{\theta}_{c}$:
\beq
        P \left( \bar{\theta}_{T^{\delta}_{b}} < c \right) &=& P \left( T^{\delta}_{b} < T^{\theta}_{c} \right) \ \ = \ \ P \left( b < \bar{\delta}_{T^{\theta}_{c}} \right) \nonumber \\
                    &=& P \left( b < \sqrt{T^{\theta}_{c}} |N| \right) \nonumber \\
                    &=& P \left( \frac{b}{\sqrt{T^{\theta}_{c}}} < |N| \right) \nonumber \\
                    &=& \sqrt{\frac{2}{\pi}} E\left[ \int^{\infty}_{b/\sqrt{T^{\theta}_{c}}} dy \: e^{-y^{2}/2} \right], \label{iiia}
\eeq
while, on the RHS of (iii):
\beq
        P \left( |C_{a(b)}| < c \right) = 2\int^{c}_{0} \frac{a(b) \; dy}{\pi(a^{2}(b)+y^{2})}\stackrel{y=a(b)h}{=} \frac{2}{\pi}\int^{c/a(b)}_{0} \frac{dh}{1+h^{2}} \label{iiib}.
\eeq
Taking derivatives in (\ref{iiia}) and (\ref{iiib}) with respect to $b$ and changing the variables
$b=\sqrt{x}$, we obtain Proposition \ref{laplace}.
\end{rem}
%%%%%%%%%%%%%%%%%%%%%%%%%%%%%%%%%%%%%%%%%%%%%%%%%%%%%%%%%%%%%%%%

%%%%%%%%%%%%%%%%%%%%%%%%%%%%%%%%%%%%%%%%%%%%%%%%%%%%%%%%%%%%%%%%
\subsection{Recovering Spitzer's theorem}
%%%%%%%%%%%%%%%%%%%%%%%%%%%%%%%%%%%%%%%%%%%%%%%%%%%%%%%%%%%%%%%%
The identity (ii) in Proposition \ref{law} is reminiscent of
Williams' remark (see \cite{Wil74,MeY82}), that:
\beq\label{williams}
    H_{T^{R}_{r}} \stackrel{(law)}{=} T^{\delta}_{\log r},
\eeq
where here $R$ starts from $1$ and $\delta$ starts from $0$
(in fact, this is a consequence of (\ref{skew-product2}) ). For a
number of variants of (\ref{williams}), see \cite{Yor85,MaY08}. This was D. Williams' starting point for
a non-computational proof of Spitzer's result (\ref{spi2}). We
note that in (ii), $T^{\delta}_{b}$ is independent of the process
$(\theta_{u},u\geq0)$ while in (\ref{williams}) $T^{R}_{r}$
depends on $(\theta_{u},u\geq0)$. Actually, we can mimic Williams'
"pinching method" to derive Spitzer's theorem (\ref{spi2}) from
(ii) in Proposition \ref{law}.
%%%%%%%%%%%%%%%%%%%%%%%%%%%%%%%%%%%%%%%%%%%%%%%%%%%%%%%%%%%%%%%%
\begin{prop}\label{spitzer} \textbf{(A new proof of Spitzer's theorem)}\\
As $t\rightarrow \infty$, $\theta_{T^{\delta}_{\sqrt{t}}} -
\theta_{t}$ converges in law, which implies that:
\beq\label{thetadifference}
    \frac{1}{\log t} \left( \theta_{T^{\delta}_{\sqrt{t}}} - \theta_{t} \right) \overset{{(P)}}{\underset{t\rightarrow\infty}\longrightarrow} 0,
\eeq
and, in turn, implies Spitzer's theorem (see formula
(\ref{spi2}) ):
\beqq
\frac{2}{\log t} \; \theta_{t}
\overset{{(law)}}{\underset{t\rightarrow\infty}\longrightarrow} C_{1}.
\eeqq
\end{prop}
%%%%%%%%%%%%%%%%%%%%%%%%%%%%%%%%%%%%%%%%%%%%%%%%%%%%%%%%%%%%%%%%
%%%%%%%%%%%%%%%%%%%%%%%%%%%%%%%%%%%%%%%%%%%%%%%%%%%%%%%%%%%%%%%%
\textbf{Proof of Proposition \ref{spitzer}} From equation (ii) of
Proposition \ref{law} we note:
\beqq
\frac{1}{\log b} \; \theta_{T^{\delta}_{b}} \stackrel{(law)}{=} \frac{C_{a(b)}}{\log
b} \overset{{(law)}}{\underset{b\rightarrow\infty}\longrightarrow} C_{1}.
\eeqq
So, for $b=\sqrt{t}$ we have:
\beqq
\frac{2}{\log t} \; \theta_{T^{\delta}_{\sqrt{t}}}
\overset{{(law)}}{\underset{b\rightarrow\infty}\longrightarrow} C_{1}.
\eeqq
On the other hand, following Williams' "pinching
method", we note that:
\beqq
\frac{1}{\log t} \left(
\theta_{T^{\delta}_{\sqrt{t}}} - \theta_{t} \right)
\overset{{(law)}}{\underset{t\rightarrow\infty}\longrightarrow} 0,
\eeqq
since $Z_{u}=x_{0}+Z^{(0)}_{u}$ and also, as we change
variables $u=tv$ and we use the scaling property, we obtain:
\beqq
\theta_{T^{\delta}_{\sqrt{t}}} - \theta_{t} \equiv \mathrm{Im}
\left( \int^{T^{\delta}_{\sqrt{t}}}_{t} \frac{dZ_{u}}{Z_{u}}
\right) \overset{{(law)}}{\underset{t\rightarrow\infty}\longrightarrow} \mathrm{Im}
\left(\int^{T^{\delta}_{1}}_{1} \frac{dZ^{(0)}_{v}}{Z^{(0)}_{v}}
\right).
\eeqq
Here, the limit variable is -in our opinion- of no
other interest than its existence which implies
(\ref{thetadifference}), hence (\ref{spi2}).
\begin{flushright}
  $\Box$
\end{flushright}
%%%%%%%%%%%%%%%%%%%%%%%%%%%%%%%%%%%%%%%%%%%%%%%%%%%%%%%%%%%%%%%%

%%%%%%%%%%%%%%%%%%%%%%%%%%%%%%%%%%%%%%%%%%%%%%%%%%%%%%%%%%%%%%%%
\subsection{On the distributions of $T^{\theta}_{c}\equiv T^{\theta}_{-\infty,c}$ and $T^{\theta}_{-c,c}$}
%%%%%%%%%%%%%%%%%%%%%%%%%%%%%%%%%%%%%%%%%%%%%%%%%%%%%%%%%%%%%%%%
%%%%%%%%%%%%%%%%%%%%%%%%%%%%%%%%%%%%%%%%%%%%%%%%%%%%%%%%%%%%%%%%
\begin{prop}\label{Tc} The asymptotic equivalence:
\beq\label{Tthetac}
  \left(\log t\right) \; P(T^{\theta}_{c}>t) \overset{{t\rightarrow \infty}}{\longrightarrow} (4c)/\pi \ ,
\eeq
holds. \\
As a consequence, for $\eta>0$, $E[(\log
T^{\theta}_{c})_{+}^{\eta}] < \infty$ if and only if $\eta <1$
(where $(\cdot)_{+}$ denotes the positive part).
\end{prop}
%%%%%%%%%%%%%%%%%%%%%%%%%%%%%%%%%%%%%%%%%%%%%%%%%%%%%%%%%%%%%%%%
%%%%%%%%%%%%%%%%%%%%%%%%%%%%%%%%%%%%%%%%%%%%%%%%%%%%%%%%%%%%%%%%
\textbf{Proof of Proposition \ref{Tc}} $\left.a\right)$ We rely upon the
asymptotic distribution of $H_{t} \equiv
\int^{t}_{0}\frac{ds}{\left|Z_{s}\right|^{2}}$ which is given by
\cite{ReY99}:
\beq\label{Hasymptotic}
  \frac{4 H_{t}}{(\log t)^{2}}  \overset{{(law)}}{\underset{t\rightarrow\infty}\longrightarrow} T_{1}^{\beta} \equiv \inf \{ t: \beta_{t}=1 \} ,
\eeq
or equivalently:
\beq\label{Hasymptotic2}
  \frac{\log t}{2 \sqrt{H_{t}}}  \overset{{(law)}}{\underset{t\rightarrow\infty}\longrightarrow} \left| N \right|,
\eeq
where $N$ is a standard Gaussian variable $\mathcal{N}(0,1)$.\\
We note that, from the representation (\ref{skew-product2}) of
$\theta_{t}$, the result (\ref{Hasymptotic}) is equivalent to
Spitzer's theorem \cite{Spi58}:
\beq\label{spi3}
  \frac{2 \theta_{t}}{\log t}  \overset{{(law)}}{\underset{t\rightarrow\infty}\longrightarrow} C_{1} \stackrel{(law)}{=} \gamma_{T_{1}^{\beta}},
\eeq
where $C_{1}$ is a standard Cauchy variable.\\
$\left.b\right)$ We shall now use this, in order to deduce Proposition \ref{Tc}. We
denote $S_{t}^{\theta} \equiv \sup_{s \leq t} \theta_{s} \equiv
S_{H_{t}}^{\gamma}$ and we note that (from scaling):
\beq\label{Tthetac2}
    P(T^{\theta}_{c} \geq t) =  P(S^{\gamma}_{H_{t}} \leq c) = P(\sqrt{H_{t}}S^{\gamma}_{1} \leq c),
\eeq
since $\gamma$ and $H$ are independent. Thus, we have (since
$S^{\gamma}_{1} \stackrel{(law)}{=} | N | $ and by making the
change of variable $x = \frac{cy}{\sqrt{H_{t}}}$):
\beq\label{Tthetac3}
    P(T^{\theta}_{c} \geq t) &=&  \sqrt{\frac{2}{\pi}} E \left[ \int_{0}^{c/\sqrt{H_{t}}} dx \; e^{-\frac{x^{2}}{2}} \right] \nonumber \\               &=&  \sqrt{ \frac{2}{\pi}} \ c \  E \left[ \int_{0}^{1}
\frac{dy}{\sqrt{H_{t}}} \exp \left( -\frac{c^{2} y^{2}}{2 H_{t}}
\right) \right].
\eeq
Thus, we now deduce from
(\ref{Hasymptotic2}) that:
\beq\label{Tthetac4}
    \frac{\log t}{2}  P(T^{\theta}_{c} \geq t) \overset{{t\rightarrow \infty}}{\longrightarrow}  \sqrt{ \frac{2}{\pi}} \; c \; E \left[ | N | \right] = \frac{2}{\pi} \; c .
\eeq
which is precisely (\ref{Tthetac}). \\
It is now elementary to deduce from (\ref{Tthetac4}) that: for
$\eta>0$:
\beqq
E[(\log T^{\theta}_{c})_{+}^{\eta}] < \infty \
\Leftrightarrow \ 0< \eta < 1 ,
\eeqq
since (\ref{Tthetac4}) is equivalent to:
\beq\label{logTthetac}
    u \; P(\log T^{\theta}_{c} > u) \overset{{u\rightarrow \infty}}{\longrightarrow} \left( \frac{4c}{\pi} \right).
\eeq
Consequently, Fubini's theorem  yields:
\beqq
E \left[(\log T^{\theta}_{c})_{+}^{\eta} \right] = \int^{\infty}_{0} du \;  \eta \; u^{\eta-1} \  P(\log T^{\theta}_{c} > u),
\eeqq
and from (\ref{logTthetac}) this is finite if and only if:
\beqq
\int^{\infty}_{\cdot} du \ u^{\eta-2} < \infty
\Leftrightarrow \eta < 1.
\eeqq
So, $ E[(\log T^{\theta}_{c})_{+}^{\eta}]< \infty \ \Leftrightarrow \ 0< \eta <1$.
\begin{flushright}
  $\Box$
\end{flushright}
%%%%%%%%%%%%%%%%%%%%%%%%%%%%%%%%%%%%%%%%%%%%%%%%%%%%%%%%%%%%%%%%
Now we give several examples of random times $T:C(\mathbb{R}_{+},\mathbb{R}) \rightarrow
\mathbb{R}_{+}$ which may be studied quite similarly to $T^{\theta}_{c}$. \\
For such times $T$, it will always be true that:
$H_{T(\theta)}=T(\gamma)$ is equivalent to $T(\theta)=A_{T(\gamma)}$,
defined with respect to $Z$, issued from $x_{0}\neq0$.
Using Bougerol's identity, we obtain:
\beq\label{bougerolgeneral}
    \sinh(\beta_{T(\gamma)}) \stackrel{(law)}{=} \hat{\beta}_{A_{T(\gamma)}} = \hat{\beta}_{(T(\theta))}.
\eeq
where $(\hat{\beta}_{u},u\geq0)$ is a 1-dimensional Brownian
motion independent of $(\beta,\gamma)$ (or equivalently, of $Z$).
Consequently, denoting by $h_{T}$ the density of
$\beta_{T(\gamma)}$, we deduce from (\ref{bougerolgeneral}) that:
\beq\label{laplacetransform3}
E \left[ \frac{1}{\sqrt{2\pi
T(\theta)}} \exp ( -\frac{x^{2}}{2T(\theta)} ) \right] =
\frac{1}{\sqrt{1+x^{2}}} \; h_{T}(\log(x+\sqrt{1+x^{2}})),
\eeq
or equivalently, changing $x$ in $\sqrt{x}$, we obtain:
\beq\label{laplacetransform4}
E \left[ \frac{1}{\sqrt{2\pi
T(\theta)}} \exp ( -\frac{x}{2T(\theta)} ) \right] =
\frac{1}{\sqrt{1+x}} \; h_{T}(\log(\sqrt{x}+\sqrt{1+x})).
\eeq
In a number of cases, $h_{T}$ is known explicitly, for example: \\
(i)
\beqq
T(\gamma)=T^{\gamma}_{-d,c} \Leftrightarrow
T(\theta)=\int^{T^{\gamma}_{-d,c}}_{0} ds \; \exp\left(2\beta_{s}
\right) =T^{\theta}_{-d,c}.
\eeqq
Thus:
\beq\label{laplacetransform5}
E \left[ \frac{1}{\sqrt{2\pi T^{\theta}_{-d,c}}} \exp (
-\frac{x}{2T^{\theta}_{-d,c}} ) \right] = \frac{1}{\sqrt{1+x}} \;
h_{-d,c}(\log(\sqrt{x}+\sqrt{1+x})),
\eeq
where $h_{-d,c}$ is the
density of the variable $\beta_{T^{\gamma}_{-d,c}}$. The law of
$\beta_{T^{\gamma}_{-d,c}}$ may be obtained from its
characteristic function which is given by
\cite{ReY99}, page 73:
\beqq
E \left[ \exp(i\lambda \beta_{T^{\gamma}_{-d,c}}) \right] &=& E \left[ \exp(-\frac{\lambda^{2}}{2} T^{\gamma}_{-d,c}) \right] \\
&=& \frac{\cosh(\frac{\lambda}{2}(c-d))}{\cosh(\frac{\lambda}{2}(c+d))}.
\eeqq
In particular, for $c=d$, we recover the very classical
formula:
\beqq
    E \left[ \exp(i\lambda \beta_{T^{\gamma}_{-c,c}}) \right] = \frac{1}{\cosh(\lambda c)}.
\eeqq
It is well known that \cite{Lev80,BiY87}:
\beq\label{Fourier}
E \left[ \exp( i\lambda \beta_{T^{\gamma}_{-c,c}}) ) \right] &=& \frac{1}{\cosh(\lambda c)} = \frac{1}{\cosh(\pi \lambda \frac{c}{\pi})} \nonumber \\
&=& \int^{\infty}_{-\infty} e^{i \left( \frac{\lambda c}{\pi} \right) x  } \frac{1}{2\pi} \frac{1}{\cosh(\frac{x}{2})} dx \nonumber \\
&\overset{y=\frac{cx}{\pi}}{=}& \int^{\infty}_{-\infty} e^{i\lambda y}  \frac{1}{2\pi} \frac{\frac{\pi}{c}}{\cosh(\frac{y \pi}{2c})} dy \nonumber \\
&=& \int^{\infty}_{-\infty} e^{i \lambda y} \frac{1}{2c}
\frac{1}{\cosh(\frac{y \pi}{2c})} dy.
\eeq
Hence, the density of
$\beta_{T^{\gamma}_{-c,c}}$ is:
\beqq
    h_{-c,c}(x)= \left( \frac{1}{2c} \right) \frac{1}{\cosh(\frac{x\pi}{2c})} = \left( \frac{1}{c} \right) \frac{1}{e^{\frac{x\pi}{2c}} + e^{-\frac{x\pi}{2c}}},
\eeqq
and
\beqq
    h_{-c,c} \left( \log (\sqrt{x}+\sqrt{1+x}) \right)= \left( \frac{1}{c} \right) \frac{1}{(\sqrt{x}+\sqrt{1+x})^{\zeta} + (\sqrt{x}+\sqrt{1+x})^{-\zeta}},
\eeqq
where $\zeta=\frac{\pi}{2c}$. However using:
\beq\label{pho}
    (\sqrt{x}+\sqrt{1+x})^{-\zeta} = (\sqrt{1+x}-\sqrt{x})^{\zeta},
\eeq
we obtain:
\beq\label{h2}
  h_{-c,c} \left( \log (\sqrt{x}+\sqrt{1+x}) \right)= \left( \frac{1}{c} \right) \frac{1}{(\sqrt{x}+\sqrt{1+x})^{\zeta} + (\sqrt{1+x} - \sqrt{x})^{\zeta}}. \nonumber \\
\eeq
So we deduce that (for $c=d$):
\beq\label{laplacetransform6}
&& E \left[ \frac{1}{\sqrt{2\pi T^{\theta}_{-c,c}}} \exp (
-\frac{x}{2T^{\theta}_{-c,c}} ) \right] \nonumber \\
 &=& \left( \frac{1}{c} \right) \left(
\frac{1}{\sqrt{1+x}} \right)
\frac{1}{(\sqrt{x}+\sqrt{1+x})^{\zeta}+(\sqrt{1+x}-\sqrt{x})^{\zeta}}.
\eeq
\begin{flushright}
  $\Box$
\end{flushright}
%%%%%%%%%%%%%%%%%%%%%%%%%%%%%%%%%%%%%%%%%%%%%%%%%%%%%%%%%%%%%%%%
(ii) As a second example of a random time $T$, let us consider the
time introduced in \cite{Val92}, \cite{ChY03}, exercise 6.2, p.
178 (we use a slightly different notation). Let
$(\beta_{t},t\geq0)$ be a real valued Brownian motion and define,
for $c>0$:
\beqq
    T(\theta)&\equiv& T^{\hat{\theta}}_{c}=\inf \left\{ t: \underset{s\leq t}{\sup} \ \theta_{s}-\underset{s\leq t}{\inf} \ \theta_{s}=c \right\}, \\
    T(\gamma)&\equiv& T^{\hat{\gamma}}_{c}=\inf \left\{ t: \underset{s\leq t}{\sup} \ \gamma_{s}-\underset{s\leq t}{\inf} \ \gamma_{s}=c \right\} .
\eeqq
Thus, from the skew-product representation (\ref{skew-product}), $\theta_{u}\equiv\gamma_{H_{u}}$, by
replacing $u=T^{\hat{\theta}}_{c}$, we obtain:
\beqq
H_{T^{\hat{\theta}}_{c}}=T^{\hat{\gamma}}_{c} \Rightarrow T^{\hat{\theta}}_{c}=\int^{T^{\hat{\gamma}}_{c}}_{0} ds \;
\exp\left(2\beta_{s} \right)\equiv A_{T^{\hat{\gamma}}_{c}}.
\eeqq
Thus:
\beq\label{laplacetransformex2}
    E \left[ \frac{1}{\sqrt{2\pi T^{\hat{\theta}}_{c}}} \exp ( -\frac{x}{2T^{\hat{\theta}}_{c}} ) \right] = \frac{1}{\sqrt{1+x}} \; h_{c}(\log(\sqrt{x}+\sqrt{1+x})),
\eeq
where $h_{c}$ is the density of the variable
$\beta_{T^{\hat{\gamma}}_{c}}$. The law of $\beta_{T^{\hat{\gamma}}_{c}}$ may be
obtained from its characteristic function which is given by \cite{BiY87,ChY03}:
\beq\label{Fourierex2}
    E \left[ \exp( i\lambda \beta_{T^{\hat{\gamma}}_{c}} ) \right] &=& E \left[ \exp( -\frac{\lambda^{2}}{2} T^{\hat{\gamma}}_{c} ) \right] = \frac{1}{(\cosh(\lambda \frac{c}{2}))^{2}} = \frac{1}{\left(\cosh(\pi \lambda \frac{c}{2\pi})\right)^{2}} \nonumber \\
    &=& \int^{\infty}_{-\infty} e^{i \left( \frac{\lambda c}{2\pi} \right) x  } \frac{1}{2\pi} \frac{x}{\sinh(\frac{x}{2})} dx \nonumber \\
    &\overset{y=\frac{cx}{2\pi}}{=}& \int^{\infty}_{-\infty} e^{i\lambda y}  \frac{1}{2\pi} \frac{\frac{2\pi y}{c}}{\sinh(\frac{\pi y}{c})} \frac{2\pi}{c} dy \nonumber \\
    &=& \int^{\infty}_{-\infty} e^{i \lambda y} \frac{2\pi}{c^{2}} \frac{y}{\sinh(\frac{\pi y}{c})} dy.
\eeq
So, the density of $\beta_{T^{\hat{\gamma}}_{c}}$ is:
\beqq
h_{c}(y)= \left( \frac{2\pi y}{c^{2}} \right) \frac{1}{\sinh(\frac{\pi
y}{c})} = \frac{4\pi}{c^{2}} \frac{y}{e^{\frac{\pi y}{c}} -
e^{-\frac{\pi y}{c}}},
\eeqq
and
\beqq
h_{c} \left( \log (\sqrt{x}+\sqrt{1+x}) \right)= \frac{4\pi}{c^{2}}
\frac{\log (\sqrt{x}+\sqrt{1+x})}{(\sqrt{x}+\sqrt{1+x})^{\hat{\zeta}} -
(\sqrt{x}+\sqrt{1+x})^{-\hat{\zeta}}},
\eeqq
where $\hat{\zeta}=\frac{\pi}{c}$. Thus:
\beq\label{laplacetransformfinex2}
&& E \left[ \frac{1}{\sqrt{2\pi T^{\hat{\theta}}_{c}}} \exp ( -\frac{x}{2T^{\hat{\theta}}_{c}} ) \right] \nonumber \\
&& = \frac{4\pi}{c^{2}} \frac{1}{\sqrt{1+x}}
\frac{\log (\sqrt{x}+\sqrt{1+x})}{(\sqrt{x}+\sqrt{1+x})^{\hat{\zeta}}-(\sqrt{1+x}-\sqrt{x})^{\hat{\zeta}}}.
\eeq
We note that this study may be related to \cite{PiY03};
and more precisely $\beta_{T^{\hat{\gamma}}_{c}}$ and
$T^{\hat{\gamma}}_{c}$ correspond to the variables $C_{2}$ and
$\hat{C_{2}}$ respectively (see e.g. Table 6 in p. 312).
\begin{flushright}
  $\Box$
\end{flushright}
%%%%%%%%%%%%%%%%%%%%%%%%%%%%%%%%%%%%%%%%%%%%%%%%%%%%%%%%%%%%%%%%
Let us now return to the case of $T^{\theta}_{-c,c}$ (example
(i)). More specifically, we shall obtain its density function
$f(t)$.
%%%%%%%%%%%%%%%%%%%%%%%%%%%%%%%%%%%%%%%%%%%%%%%%%%%%%%%%%%%%%%%%
\begin{prop}\label{density}
The density function $f$ of $T^{\theta}_{-c,c}$ is given by:
\beq\label{densityfunction}
  f(t)=\frac{1}{\sqrt{2}c} \sum ^{\infty}_{k=0} (-1)^{k} \frac{\Gamma (\nu_{k})}{\Gamma(2\nu_{k})} \frac{1}{\sqrt{t}} e^{-\frac{1}{4t}} M_{\frac{1}{2},\nu_{k}}(\frac{1}{2t}),
\eeq
where $M_{a,b}(\cdot)$ is the Whittaker function with
parameters $a,b$. Equivalently:
\beq\label{densityfunction2}
f(t) = \frac{\sqrt{2}}{c} \sum ^{\infty}_{k=0} (-1)^{k}
\frac{1}{\sqrt{t}} e^{-\frac{1}{2t}} \left( \frac{1}{2t}
\right)^{\nu_{k} + \frac{1}{2}} \nu_{k} \sum ^{\infty}_{n=0}
\frac{\Gamma (\nu_{k}+n)}{\Gamma(2 \nu_{k}+n+1)} \frac{1}{n!}
\left( \frac{1}{2t} \right) ^{n}, \nonumber \\
\eeq
where $\nu_{k} = \frac{\pi}{4c}(2k +1)$.
\end{prop}
%%%%%%%%%%%%%%%%%%%%%%%%%%%%%%%%%%%%%%%%%%%%%%%%%%%%%%%%%%%%%%%%
%%%%%%%%%%%%%%%%%%%%%%%%%%%%%%%%%%%%%%%%%%%%%%%%%%%%%%%%%%%%%%%%
\textbf{Proof of Proposition \ref{density}}
The following calculation relies upon a private note by A. Comtet
\cite{Com06}.We denote:
\beqq
\varphi_{\zeta}(x) =
(\sqrt{x}+\sqrt{1+x})^{\zeta}+(\sqrt{1+x}-\sqrt{x})^{\zeta}.
\eeqq
Noting:
\beq\label{cosh}
\sqrt{1+x} = \cosh \frac{y}{2}
\Longleftrightarrow y = 2 \arg\cosh(\sqrt{1+x}),
\eeq
we get:
\beqq
\varphi_{\zeta}(x) &=& (\sinh \frac{y}{2} + \cosh \frac{y}{2})^{\zeta}+(\cosh \frac{y}{2}-\sinh \frac{y}{2})^{\zeta} \\
&=& 2\cosh \frac{y\zeta}{2}.
\eeqq
Thus, from
(\ref{laplacetransform6}), we have:
\beq\label{II}
II := E \left[
\frac{1}{\sqrt{2\pi T^{\theta}_{-c,c}}} \exp (
-\frac{x}{2T^{\theta}_{-c,c}} ) \right] =
\frac{1}{\psi}\frac{1}{\cosh \frac{y}{2}} \frac{1}{\cosh \frac{\pi y}{2\psi}},
\eeq
where $\psi=2c$. However, expanding $\cosh
\frac{\pi y}{2\psi}$, we get:
\beqq
\frac{1}{\cosh \frac{\pi
y}{2\psi}} = 2 \frac{e^{-\frac{\pi y}{2\psi}}}{1+e^{-\frac{\pi
y}{\psi}}} = 2 \sum^{\infty}_{k=0} \left( -e^{-\frac{\pi y}{\psi}}
\right)^{k} e^{-\frac{\pi y}{2\psi}},
\eeqq
and from (\ref{II}),
we deduce that:
\beqq
II &=& \sum^{\infty}_{k=0} \frac{2}{\psi} \frac{(-1)^{k}}{\cosh \frac{y}{2}} e^{-\frac{\pi}{2\psi}(2k +1)y} \\
&=& \sum^{\infty}_{k=0} \frac{4(-1)^{k}}{\psi \sqrt{2} \sqrt{2\sinh \frac{y}{2} \cosh \frac{y}{2}}} \sqrt{\frac{\sinh \frac{y}{2}}{\cosh \frac{y}{2}}} e^{-\nu_{k} y} \\
&=& \sum^{\infty}_{k=0} \frac{4(-1)^{k}}{\psi \sqrt{2}
\sqrt{2\sinh \frac{y}{2} \cosh \frac{y}{2}}} \sqrt{\tanh
\frac{y}{2}} e^{-\nu_{k} y},
\eeqq
where $\nu_{k}=\frac{\pi}{2\psi}(2k +1)$. \\
From (\ref{cosh}), we have $1+x = \cosh^{2}
\frac{y}{2} \Longleftrightarrow x = \sinh^{2} \frac{y}{2}$, thus:
\beqq
    (\tanh \frac{y}{2})^{1/2} = \sqrt{\frac{\sinh \frac{y}{2}}{\cosh \frac{y}{2}}} = \left( \frac{\sqrt{x}}{\sqrt{1+x}} \right)^{1/2} = \left( \frac{x}{1+x} \right)^{1/4}.
\eeqq
Moreover, we know that
(see \cite{AbSt70}, equation 8.6.10, or \cite{Leb72}):
\beqq
i \sqrt{\frac{\pi}{2\sinh y}} e^{-\nu_{k} y} = Q^{1/2}_{\nu_{k} -
1/2}(\cosh y),
\eeqq
where $\left\{Q^{a}_{b}(\cdot) \right\}$ is
the family of Legendre functions and $\cosh y = 2x+1 $. So, we
deduce:
\beq\label{II2}
    II = \sum ^{\infty}_{k=0} \frac{4(-i)}{\psi \sqrt{\pi}} (-1)^{k} \left( \frac{x}{1+x} \right)^{1/4}  Q^{1/2}_{\nu_{k} - 1/2}(2x+1).
\eeq
By using formula 7.621.9, page 864 in \cite{GrR65}:
\beq\label{whittaker}
\int^{\infty}_{0} e^{-sw}
M_{l,\nu_{k}}(w) \frac{dw}{w} = \frac{2\Gamma
(1+2\nu_{k}) \; e^{-i\pi l}}{\Gamma(\frac{1}{2}+\nu_{k}+l)} \left(
\frac{s-\frac{1}{2}}{s+\frac{1}{2}} \right)^{l/2} Q^{l}_{\nu_{k} -
1/2}(2s),
\eeq
with: $l=\frac{1}{2}$, $\nu_{k}=\frac{\pi}{2\psi}(2k
+1)$, $s=x+\frac{1}{2}$ and $M_{\cdot,\cdot}(\cdot)$ denoting the
Whittaker function, which is defined as:
\beqq
    M_{a,b}(w) = w^{b + \frac{1}{2}} e^{-\frac{1}{2}w} \frac{\Gamma (2b+1)}{\Gamma(\frac{1}{2}+b-a)} \sum ^{\infty}_{n=0} \frac{\Gamma (\frac{1}{2}+b-a+n)}{\Gamma(2b+1+n)} \frac{w^{n}}{n!}.
\eeqq
we have:
\beq\label{whittaker2}
-2i \frac{\Gamma (1+2\nu_{k})}{\Gamma(1+\nu_{k})} \left( \frac{x}{1+x}
\right)^{1/4}  Q^{1/2}_{\nu_{k} - 1/2}(2x+1) = \int^{\infty}_{0}
e^{-sw} M_{1/2,\nu_{k}}(w) \frac{dw}{w}. \nonumber \\
\eeq
From (\ref{II2}) and
by changing the variable $w = \frac{1}{2t}$, we deduce:
\beq\label{whittaker3}
II &=& \sum ^{\infty}_{k=0} \frac{2}{\psi \sqrt{\pi}} (-1)^{k} \frac{\Gamma (\nu_{k}+1)}{\Gamma(2\nu_{k}+1)} \int^{\infty}_{0} \frac{dw}{w} \exp \left(-w(x+\frac{1}{2}) \right) M_{1/2,\nu_{k}}(w) \nonumber \\
&=& \sum ^{\infty}_{k=0} \int^{\infty}_{0} \frac{dt}{t} \frac{2}{\psi \sqrt{\pi}} (-1)^{k} \frac{\Gamma (\nu_{k}+1)}{\Gamma(2\nu_{k}+1)} \exp \left(-\frac{1}{4t}-\frac{x}{2t} \right) M_{1/2,\nu_{k}}(\frac{1}{2t}) . \nonumber \\
\eeq
By using the equations (\ref{II}) and (\ref{whittaker3}), we
conclude:
\beq\label{whittaker4}
 &&   E \left[ \frac{1}{\sqrt{2\pi T^{\theta}_{-c,c}}} \exp \left( -\frac{x}{2T^{\theta}_{-c,c}} \right) \right] \nonumber \\
 &=& \sum ^{\infty}_{k=0} \int^{\infty}_{0} \frac{dt}{t} \frac{2}{\psi \sqrt{\pi}} (-1)^{k} \frac{\Gamma (\nu_{k}+1)}{\Gamma(2\nu_{k}+1)} \exp \left(-\frac{1}{4t}-\frac{x}{2t} \right) M_{\frac{1}{2},\nu_{k}}(\frac{1}{2t}) \nonumber \\
 &=& \sum ^{\infty}_{k=0}
\int^{\infty}_{0} \frac{dt}{t} \frac{2}{\psi \sqrt{\pi}} (-1)^{k}
\frac{\Gamma (\frac{\pi}{4c}(2k +1)+1)}{\Gamma(2\frac{\pi}{4c}(2k
+1)+1)} \exp \left(-\frac{1}{4t}-\frac{x}{2t} \right)
M_{\frac{1}{2},\frac{\pi}{4c}(2k +1)}(\frac{1}{2t}). \nonumber \\
\eeq
Thus, the density function $f$ of $T^{\theta}_{-c,c}$ is given by:
\beq
        f(t) &=&  \frac{2 \sqrt{2}}{\psi} \sum ^{\infty}_{k=0} (-1)^{k} \frac{\Gamma (\nu_{k}+1)}{\Gamma(2\nu_{k}+1)} \frac{1}{\sqrt{t}} e^{-\frac{1}{4t}} M_{\frac{1}{2},\nu_{k}}(\frac{1}{2t}) \label{densityfunction3} \\
                 &=& \frac{\sqrt{2}}{c} \sum ^{\infty}_{k=0} (-1)^{k} \frac{\Gamma (\frac{\pi}{4a}(2k +1)+1)}{\Gamma(\frac{\pi}{2a}(2k +1)+1)} \frac{1}{\sqrt{t}} e^{-\frac{1}{4t}} M_{\frac{1}{2},\frac{\pi}{4a}(2k +1)}(\frac{1}{2t}) \label{densityfunction4} \\
                 &=& \frac{\sqrt{2}}{c} \sum ^{\infty}_{k=0} (-1)^{k} \frac{\nu_{k} \Gamma (\nu_{k})}{2\nu_{k} \Gamma(2\nu_{k})} \frac{1}{\sqrt{t}} e^{-\frac{1}{4t}} M_{\frac{1}{2},\nu_{k}}(\frac{1}{2t}), \label{densityfunction5}
\eeq
where the Whittaker function
$M_{\frac{1}{2},\nu_{k}}(\frac{1}{2t})$ is:
\beq
  &&  M_{\frac{1}{2},\frac{\pi}{4c}(2k +1)}(\frac{1}{2t}) \nonumber \\
        &=& \left( \frac{1}{2t} \right)^{\frac{\pi}{4c}(2k +1) + \frac{1}{2}} e^{-\frac{1}{4t}} \frac{\Gamma (\frac{\pi}{2c}(2k +1)+1)}{\Gamma(\frac{\pi}{4c}(2k +1))} \sum ^{\infty}_{n=0} \frac{\Gamma (\frac{\pi}{4c}(2k +1)+n)}{\Gamma(\frac{\pi}{2c}(2k +1)+1+n)} \frac{1}{n!} \left( \frac{1}{2t} \right) ^{n} \nonumber \\
        &=& \left( \frac{1}{2t} \right)^{\nu_{k} + \frac{1}{2}} e^{-\frac{1}{4t}} \frac{\Gamma (2\nu_{k}+1)}{\Gamma(\nu_{k})} \sum ^{\infty}_{n=0} \frac{\Gamma (\nu_{k}+n)}{\Gamma(2 \nu_{k}+1+n)} \frac{1}{n!} \left( \frac{1}{2t} \right) ^{n} \nonumber \\
        &=& \left( \frac{1}{2t} \right)^{\nu_{k} + \frac{1}{2}} e^{-\frac{1}{4t}} (2\nu_{k}) \frac{\Gamma (2\nu_{k})}{\Gamma(\nu_{k})} \sum ^{\infty}_{n=0} \frac{\Gamma (\nu_{k}+n)}{(2\nu_{k}+n) \Gamma(2 \nu_{k}+n)} \frac{1}{n!} \left( \frac{1}{2t} \right) ^{n}. \label{whittaker5}
\eeq
Thus, from (\ref{densityfunction5}) and (\ref{whittaker5}), we deduce (\ref{densityfunction2}).\\
\begin{flushright}
  $\Box$
\end{flushright}
%%%%%%%%%%%%%%%%%%%%%%%%%%%%%%%%%%%%%%%%%%%%%%%%%%%%%%%%%%%%%%%%

Next, we present the graphs of different approximations $f_{K,N}(t)$ of $f(t)$, in (\ref{densityfunction2}),
where $f_{K,N}$ denotes the sum in the series in (\ref{densityfunction2}) of the terms for
$k\leq K$, and $n\leq N$.

%%%%%%%%%%%%%%%%%%%%%%%%%%%%%%%%%%%%%%%%%%%%%%%%%%%%%%%%%%%%%%%%
\begin{figure}
\includegraphics[width=0.75\textwidth]{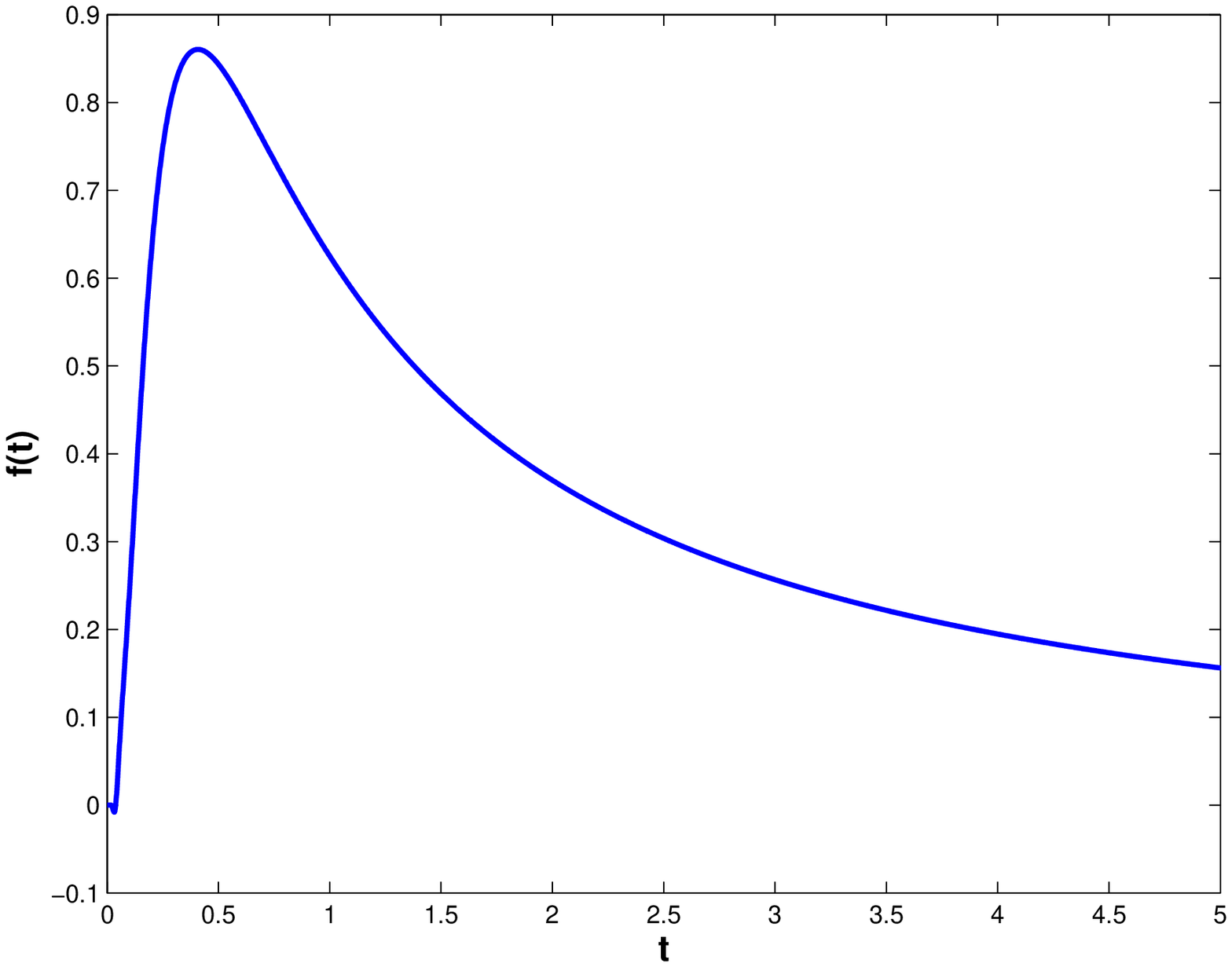}
\caption{Graph of $f_{9,9}(t)$, with $c=2\pi$. }\label{figf}
\end{figure}
%%%%%%%%%%%%%%%%%%%%%%%%%%%%%%%%%%%%%%%%%%%%%%%%%%%%%%%%%%%%%%%%

%%%%%%%%%%%%%%%%%%%%%%%%%%%%%%%%%%%%%%%%%%%%%%%%%%%%%%%%%%%%%%%%
\begin{figure}
\includegraphics[width=0.75\textwidth]{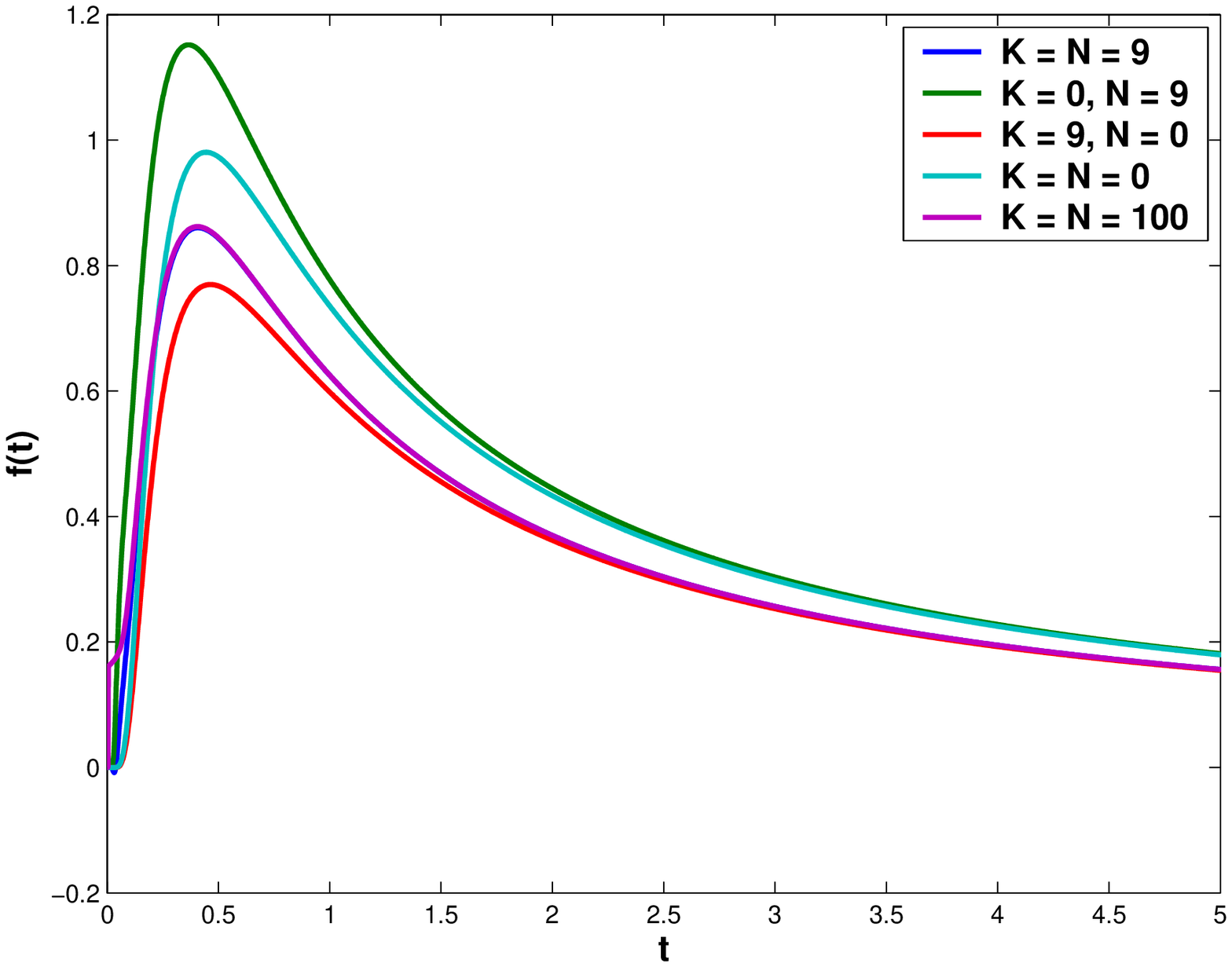}
\caption{Graph of $f_{K,N}(t)$ for several
values of $K$ and $N$, with $c=2\pi$.}\label{figf2}
\end{figure}
%%%%%%%%%%%%%%%%%%%%%%%%%%%%%%%%%%%%%%%%%%%%%%%%%%%%%%%%%%%%%%%%
\begin{rem}\label{diagram}
\begin{itemize}
    \item { Figure \ref{figf} represents the approximation of the density function $f$ with respect to the time $t$ (for $K$ and $N\leq 9$), with $c=2\pi$, whereas Figure \ref{figf2} represents the approximation of $f$ with respect to the time $t$ for several values of $k$ and $n$, with $c=2\pi$. }
    \item { From Figure \ref{figf2}, we may remark that the approximation $K$ and $N\leq 9$ is sufficiently good (comparing to the one for $K$ and $N\leq 100$). }
    \item { For the case $K$ and $N\leq 9$ it seems that locally, in a small area around 0, $f(t)<0$ which is not right. This is due to the first negative ($k=1$) term of the sum and due to the fact that we have omitted many terms. However, this is not a problem because it appears only locally. Similar irregularities have already been observed in previous articles \cite{Ish05} p.275. }
\end{itemize}
\end{rem}

%%%%%%%%%%%%%%%%%%%%%%%%%%%%%%%%%%%%%%%%%%%%%%%%%%%%%%%%%%%%%%%%
\subsection{On the first moment of $\ln\left(T^{\theta}_{-c,c}\right)$}
%%%%%%%%%%%%%%%%%%%%%%%%%%%%%%%%%%%%%%%%%%%%%%%%%%%%%%%%%%%%%%%%

This subsection is related to a result in \cite{CMY98}.
%%%%%%%%%%%%%%%%%%%%%%%%%%%%%%%%%%%%%%%%%%%%%%%%%%%%%%%%%%%%%%%%
\begin{prop}\label{logTc}
The first moment of $\ln\left(T^{\theta}_{-c,c}\right)$ has the
following integral representation:
\beq\label{logTcequation}
    E\left[\ln\left(T^{\theta}_{-c,c}\right)\right] = 2 \int^{\infty}_{0} \frac{dz}{\cosh \left(\frac{\pi z}{2}\right)} \ln\left(\sinh\left(cz\right)\right) + \ln\left(2\right) + c_{E},
\eeq
where $c_{E}=-\Gamma'(1)$ is the Euler-Mascheroni constant
(also called Euler's constant).
\end{prop}
%%%%%%%%%%%%%%%%%%%%%%%%%%%%%%%%%%%%%%%%%%%%%%%%%%%%%%%%%%%%%%%%
%%%%%%%%%%%%%%%%%%%%%%%%%%%%%%%%%%%%%%%%%%%%%%%%%%%%%%%%%%%%%%%%
\textbf{Proof of Proposition \ref{logTc}} Let us return to
equations (\ref{skew-product2}) and (\ref{TcA}). So, for
$t=T^{\theta}_{-c,c}$, we have:
\beq\label{TAln}
    \theta_{T^{\theta}_{-c,c}}=\gamma_{H_{T^{\theta}_{-c,c}}} \Longleftrightarrow H_{T^{\theta}_{-c,c}}=T^{\gamma}_{-c,c} \Longleftrightarrow T^{\theta}_{-c,c}=A_{T^{\gamma}_{-c,c}}.
\eeq
Thus, for $\varepsilon>0$:
\beqq
E\left[\left(T^{\theta}_{-c,c}\right)^{\varepsilon}\right] =
E\left[\left(A_{T^{\gamma}_{-c,c}}\right)^{\varepsilon}\right].
\eeqq
Consider $(\delta_{u},u\geq0)$ a Brownian motion,
independent of $A_{t}$. Then, Bougerol's identity and the scaling
property yield ($\mathcal{G}_{a}$ denotes a gamma variable with
parameter $a$, and $N^{2}\stackrel{(law)}{=}2 \mathcal{G}_{1/2}$):
\beqq
    E\left[\left(\sinh\left(B_{t}\right)\right)^{2\varepsilon}\right] &=& E\left[\left(\delta_{A_{t}}\right)^{2\varepsilon}\right] = E\left[A_{t}^{\varepsilon} \left(\delta_{1}\right)^{2\varepsilon}\right] \\
        &=& E\left[A_{t}^{\varepsilon}\right] \: E\left[\left(2\mathcal{G}_{1/2}\right)^{\varepsilon}\right] \\
        &=& E\left[A_{t}^{\varepsilon}\right] \: \left(2^{\varepsilon}\right) \: \frac{\Gamma\left(\frac{1}{2}+\varepsilon\right)}{\Gamma\left(\frac{1}{2}\right)},
\eeqq
because
\beqq
    E\left[\left(\mathcal{G}_{1/2}\right)^{\varepsilon}\right] &=& \int^{\infty}_{0} x^{\varepsilon+\frac{1}{2}-1} \frac{e^{-x}}{\Gamma\left(\frac{1}{2}\right)} \: dx = \frac{\Gamma\left(\frac{1}{2}+\varepsilon\right)}{\Gamma\left(\frac{1}{2}\right)}.
\eeqq
Thus, for $t=T^{\gamma}_{-c,c}$, we have:
\beq\label{shAGamma}
    E\left[\left(\sinh\left(B_{T^{\gamma}_{-c,c}}\right)\right)^{2\varepsilon}\right] = E\left[A_{T^{\gamma}_{-c,c}}^{\varepsilon}\right] \: \left(2^{\varepsilon}\right) \: \frac{\Gamma\left(\frac{1}{2}+\varepsilon\right)}{\Gamma\left(\frac{1}{2}\right)}.
\eeq
Recall that \cite{Lev80,BiY87}:
\beqq E \left[ \exp( i\lambda
B_{T^{\gamma}_{-c,c}}) ) \right] = E \left[ \exp(
-\frac{\lambda^{2}}{2} T^{\gamma}_{-c,c}) ) \right] =
\frac{1}{\cosh(\lambda c)},
\eeqq
and the density of $\beta_{T^{\gamma}_{-c,c}}$ is:
\beqq
    h_{-c,c}(y)= \left( \frac{1}{2c} \right) \frac{1}{\cosh(\frac{y\pi}{2c})} = \left( \frac{1}{c} \right) \frac{1}{e^{\frac{y\pi}{2c}} + e^{-\frac{y\pi}{2c}}}.
\eeqq
Thus, on the left hand side of (\ref{shAGamma}), we have:
\beqq
    E\left[\left(\sinh\left(B_{T^{\gamma}_{-c,c}}\right)\right)^{2\varepsilon}\right] &=& \int^{\infty}_{-\infty} \frac{dy}{2c} \:  \frac{1}{\cosh(\frac{\pi y}{2c})} \left(\sinh\left(y\right)\right)^{2\varepsilon} \\
    &=& \int^{\infty}_{0} \frac{dy}{c} \: \frac{1}{\cosh(\frac{\pi y}{2c})} \left(\sinh y\right)^{2\varepsilon}  \\
    &=& \int^{\infty}_{0} dz \: \frac{1}{\cosh(\frac{\pi z}{2})} \left(\sinh (cz)\right)^{2\varepsilon},
\eeqq
where we have made the change of variable $z=\frac{y}{c}$.
Hence, from (\ref{shAGamma}), by writing:
\beqq
E\left[A_{T^{\gamma}_{-c,c}}^{\varepsilon}\right] =
E\left[\left(T^{\theta}_{-c,c}\right)^{\varepsilon}\right] =
E\left[ e^{\varepsilon \;
\ln\left(T^{\theta}_{-c,c}\right)}\right],
\eeqq
we deduce:
\beqq
\frac{\Gamma\left(\frac{1}{2}+\varepsilon\right)}{\Gamma\left(\frac{1}{2}\right)}
E\left[ e^{\varepsilon \;
\ln\left(T^{\theta}_{-c,c}\right)}\right] =
\frac{1}{2^{\varepsilon}} \int^{\infty}_{0}
\frac{dz}{\cosh(\frac{\pi z}{2})} \left(\sinh
(cz)\right)^{2\varepsilon},
\eeqq
and by removing 1 from both sides, we obtain:
\beq\label{logepsilon}
    \frac{\Gamma\left(\frac{1}{2}+\varepsilon\right)}{\Gamma\left(\frac{1}{2}\right)} E\left[ e^{\varepsilon \; \ln\left(T^{\theta}_{-c,c}\right)}\right] - 1 = \int^{\infty}_{0} \frac{dz}{\cosh(\frac{\pi z}{2})} \left( \frac{\left( \sinh (cz)\right)^{2\varepsilon}}{2^{\varepsilon}} - 1 \right).
\eeq
On the left hand side, we apply the trivial identity
$ab-1=a(b-1)+a-1$ with
$a=\frac{\Gamma\left(\frac{1}{2}+\varepsilon\right)}{\Gamma\left(\frac{1}{2}\right)}$
and $b=E\left[ e^{\varepsilon \;
\ln\left(T^{\theta}_{-c,c}\right)}\right]$, we divide by
$\varepsilon$ and we take the limit for $\varepsilon\rightarrow0$.
Thus:
\beqq
    \frac{a(b-1)}{\varepsilon} &=& \frac{\Gamma\left(\frac{1}{2}+\varepsilon\right)}{\Gamma\left(\frac{1}{2}\right)} \frac{E\left[ e^{\varepsilon \; \ln\left(T^{\theta}_{-c,c}\right)}\right] - 1}{\varepsilon} \\
    &\stackrel{\varepsilon\rightarrow0}{\longrightarrow}& E\left[ \ln\left(T^{\theta}_{-c,c}\right)\right],
\eeqq
and:
\beqq
    \frac{a-1}{\varepsilon} &=& \frac{1}{\varepsilon} \left( \frac{\Gamma\left(\frac{1}{2}+\varepsilon\right)}{\Gamma\left(\frac{1}{2}\right)} -1 \right) = \frac{1}{\Gamma\left(\frac{1}{2}\right)} \left(\frac{\Gamma\left(\frac{1}{2}+\varepsilon\right)-\Gamma\left(\frac{1}{2}\right)}{\varepsilon} \right) \\
    &\stackrel{\varepsilon\rightarrow0}{\longrightarrow}& \frac{1}{\sqrt{\pi}} \Gamma'\left(\frac{1}{2}\right) = \frac{1}{\sqrt{\pi}} \left(-\sqrt{\pi}\right) \left(c_{E}+ 2\ln2\right) = - \left(c_{E}+ 2\ln2\right).
\eeqq
On the right hand side of (\ref{logepsilon}), we have:
\beqq
\frac{1}{\varepsilon} \left[ \left( \frac{\left( \sinh
(cz)\right)^{2}}{2}\right)^{\varepsilon} - 1 \right] =
\frac{1}{\varepsilon} \left[ \exp\left( \varepsilon \;
\ln\left(\frac{\left(\sinh (cz)\right)^{2}}{2}\right) \right) - 1
\right],
\eeqq
hence:
\beqq
&& \frac{1}{\varepsilon} \int^{\infty}_{0}
\frac{dz}{\cosh(\frac{\pi z}{2})} \left( \frac{\left( \sinh
(cz)\right)^{2\varepsilon}}{2^{\varepsilon}} - 1 \right) \\
&\stackrel{\varepsilon\rightarrow0}{\longrightarrow}&
\int^{\infty}_{0} \frac{dz}{\cosh(\frac{\pi z}{2})}
\ln\left(\frac{\left(\sinh (cz)\right)^{2}}{2}\right) \\
&=& - \ln\left(2\right)+ 2 \int^{\infty}_{0} \frac{dz}{\cosh(\frac{\pi z}{2})}
\left(\ln\left(\sinh (cz)\right)\right) \ ,
\eeqq
which finishes the proof.
\begin{flushright}
  $\Box$
\end{flushright}
%%%%%%%%%%%%%%%%%%%%%%%%%%%%%%%%%%%%%%%%%%%%%%%%%%%%%%%%%%%%%%%%

%%%%%%%%%%%%%%%%%%%%%%%%%%%%%%%%%%%%%%%%%%%%%%%%%%%%%%%%%%%%%%%%
\begin{rem}
$\left.a\right)$ We denote now:
\beq
F(c)\equiv\int^{\infty}_{0} \frac{dz}{\cosh \left(\frac{\pi z}{2}\right)} \ln\left(\sinh\left(cz\right)\right).
\eeq
Thus:
\beq
F(c)-\ln(c)\equiv \int^{\infty}_{0} \frac{dz}{\cosh \left(\frac{\pi z}{2}\right)} \ln\left(\frac{\sinh\left(cz\right)}{c}\right) \stackrel{c\rightarrow0}{\longrightarrow} \int^{\infty}_{0} \frac{dz \ \ln(z)}{\cosh \left(\frac{\pi z}{2}\right)} \approx -0.7832. \nonumber \\
\eeq
$\left.b\right)$ More generally, we denote:
\beq\label{Fcdelta}
F(c,\delta)\equiv\int^{\infty}_{0} \frac{dz}{\cosh \left(\delta z\right)} \ln\left(\sinh\left(cz\right)\right),
\eeq
and, changing the variables: $z=\frac{\pi}{2\delta}u$, we obtain:
\beq\label{FcdeltaFc}
F\left(c,\delta\right)=\left(\frac{\pi}{2\delta}\right) \int^{\infty}_{0} \frac{du}{\cosh \left(\frac{\pi}{2} u \right)} \ln\left(\sinh\left(c \frac{\pi}{2\delta} u \right)\right)= \frac{\pi}{2\delta} \ F\left(c \ \frac{\pi}{2\delta}\right).
\eeq
\end{rem}
%%%%%%%%%%%%%%%%%%%%%%%%%%%%%%%%%%%%%%%%%%%%%%%%%%%%%%%%%%%%%%%%

\newpage

%%%%%%%%%%%%%%%%%%%%%%%%%%%%%%%%%%%%%%%%%%%%%%%%%%%%%%%%%%%%%%%%
\section{The Ornstein-Uhlenbeck case}\label{sectionOU}
%%%%%%%%%%%%%%%%%%%%%%%%%%%%%%%%%%%%%%%%%%%%%%%%%%%%%%%%%%%%%%%%

%%%%%%%%%%%%%%%%%%%%%%%%%%%%%%%%%%%%%%%%%%%%%%%%%%%%%%%%%%%%%%%%
\subsection{An identity in law for Ornstein-Uhlenbeck processes, which is connected to Bougerol's identity}
%%%%%%%%%%%%%%%%%%%%%%%%%%%%%%%%%%%%%%%%%%%%%%%%%%%%%%%%%%%%%%%%
Consider the complex valued Ornstein-Uhlenbeck (OU) process:
\beq\label{OUequation}
    Z_{t} = z_{0} + \tilde{Z_{t}} - \lambda \int^{t}_{0} Z_{s} ds,
\eeq
where $\tilde{Z_{t}}$ is a complex valued Brownian motion
(BM), $z_{0}\in \mathbb{C}$ and $\lambda \geq 0$ and
$T^{(\lambda)}_{c} \equiv T^{\theta^{Z}}_{-c,c} \equiv \inf
\left\{t\geq 0 : \left|\theta^{Z}_{t}\right|=c \right\}$
($\theta^{Z}_{t}$ is the continuous winding process associated to
$Z$) denoting the first hitting time of the symmetric conic
boundary of angle $c$ for $Z$. It is well known that \cite{ReY99}:
\beq\label{OUeB}
    Z_{t} &=& e^{-\lambda t} \left( z_{0} + \int^{t}_{0} e^{\lambda s} d\tilde{Z_{s}} \right) \nonumber \\
          &=& e^{-\lambda t} \left( \mathbb{B}_{\alpha_{t}} \right),
\eeq
where, in the second equation, with the help of
Dambis-Dubins-Schwarz Theorem, $\left(\mathbb{B}_{t},t\geq
0\right)$ is a complex valued Brownian motion starting from
$z_{0}$ and
\beqq
    \alpha_{t}=\int^{t}_{0} e^{2\lambda s} ds = \frac{e^{2\lambda t}-1}{2 \lambda} \ \ .
\eeqq
We are interested in the study of the continuous winding
process
$\theta^{Z}_{t}=\mathrm{Im}(\int^{t}_{0}\frac{dZ_{s}}{Z_{s}}),t\geq0$.
By applying It\^{o}'s formula to (\ref{OUeB}), we have:
\beqq
dZ_{s} = e^{-\lambda s}(-\lambda) \mathbb{B}_{\alpha_{s}} ds +
e^{-\lambda s} d\left(\mathbb{B}_{\alpha_{s}}\right).
\eeqq
We divide by $Z_{s}$ and we obtain:
\beqq
\frac{dZ_{s}}{Z_{s}} =
-\lambda \ ds + \frac{d\mathbb{B}_{\alpha_{s}}}{\mathbb{B}_{\alpha_{s}}},
\eeqq
hence:
\beqq
\mathrm{Im} \left(\frac{dZ_{s}}{Z_{s}}\right) =
\mathrm{Im}
\left(\frac{d\mathbb{B}_{\alpha_{s}}}{\mathbb{B}_{\alpha_{s}}}\right),
\eeqq
which means that:
\beqq
\theta^{Z}_{t} = \theta^{\mathbb{B}}_{\alpha_{t}} \ \ .
\eeqq
Thus, the following holds:
%%%%%%%%%%%%%%%%%%%%%%%%%%%%%%%%%%%%%%%%%%%%%%%%%%%%%%%%%%%%%%%%
\begin{prop}\label{OU} Using the previously introduced notation, we have:
\beq\label{thetaBMOU}
    \theta^{Z}_{t} = \theta^{\mathbb{B}}_{\alpha_{t}},
\eeq
and:
\beq\label{Tchat}
    T^{(\lambda)}_{c}=\frac{1}{2\lambda}\ln \left(1+2\lambda T^{\theta^{\mathbb{B}}}_{-c,c}\right),
\eeq
where $T^{\theta^{\mathbb{B}}}_{-c,c}$ is the exit time from a cone of
angle $c$ for the complex valued BM $\mathbb{B}$.
\end{prop}
%%%%%%%%%%%%%%%%%%%%%%%%%%%%%%%%%%%%%%%%%%%%%%%%%%%%%%%%%%%%%%%%
%%%%%%%%%%%%%%%%%%%%%%%%%%%%%%%%%%%%%%%%%%%%%%%%%%%%%%%%%%%%%%%%
\textbf{Proof of Proposition \ref{OU}} We define
\beq\label{Tchat2}
    T^{(\lambda)}_{c} &\equiv& T^{\theta^{Z}}_{-c,c} \equiv \inf \left\{t\geq 0 : \left|\theta^{Z}_{t}\right|=c \right\} \nonumber \\
    &=& \inf \left\{t\geq 0 : \left|\theta^{\mathbb{B}}_{\alpha_{t}}\right|=c \right\}.
\eeq
Thus, we deduce that
$\alpha_{T^{(\lambda)}_{c}}=T^{\theta^{\mathbb{B}}}_{-c,c}\equiv
T^{\theta}_{-c,c}$. However, $T^{\theta}_{-c,c}$ (the exit time
from a cone for the BM) has already been studied in the previous
chapter and we know the explicit formula of its density function
(Proposition \ref{density}). Thus:
\beq\label{Tchat3}
    T^{(\lambda)}_{c}=\alpha^{-1}\left(T^{\theta^{\mathbb{B}}}_{-c,c}\right)=\alpha^{-1}\left(T^{\theta}_{-c,c}\right),
\eeq
where $\alpha^{-1}(t)=\frac{1}{2\lambda}\ln \left(1+2\lambda
t\right)$. Consequently:
\beqq
T^{(\lambda)}_{c}=\frac{1}{2\lambda}\ln \left(1+2\lambda
T^{\theta}_{-c,c}\right),
\eeqq
and:
\beq\label{ETchat}
    E\left[T^{(\lambda)}_{c}\right]=\frac{1}{2\lambda} E\left[\ln \left(1+2\lambda T^{\theta}_{-c,c}\right) \right],
\eeq
which finishes the proof.
\begin{flushright}
  $\Box$
\end{flushright}
%%%%%%%%%%%%%%%%%%%%%%%%%%%%%%%%%%%%%%%%%%%%%%%%%%%%%%%%%%%%%%%%
From now on, for simplicity, we shall take $z_{0}=1$ (but this is
really
no restriction, as the dependency in $z_{0}$, which is exhibited in (\ref{z0dep}), is very simple). \\
The following Proposition may be considered as an extension of the
identity in law (ii) in Proposition \ref{law}, which results from
Bougerol's identity.
%%%%%%%%%%%%%%%%%%%%%%%%%%%%%%%%%%%%%%%%%%%%%%%%%%%%%%%%%%%%%%%%
\begin{prop}\label{OUbougerol}
Consider $(Z^{\lambda}_{t},t\geq 0)$ and $(U^{\lambda}_{t},t\geq
0)$ two independent Ornstein-Uhlenbeck processes, the first one complex valued and the second one real valued,
both starting from a point different from 0, and call
$T^{(\lambda)}_{b}(U^{\lambda})= \inf \left\{t\geq 0 : e^{\lambda
t } U^{\lambda}_{t}=b \right\}$. Then, an Ornstein-Uhlenbeck
extension of identity in law (ii) in Proposition \ref{law} is the
following:
\beq\label{OUbougerolequation}
    \theta^{Z^{\lambda}}_{T^{(\lambda)}_{b}(U^{\lambda})}  \stackrel{(law)}{=} C_{a(b)},
\eeq
where $a(x)= \arg \sinh (x)$.
\end{prop}
%%%%%%%%%%%%%%%%%%%%%%%%%%%%%%%%%%%%%%%%%%%%%%%%%%%%%%%%%%%%%%%%
%%%%%%%%%%%%%%%%%%%%%%%%%%%%%%%%%%%%%%%%%%%%%%%%%%%%%%%%%%%%%%%%
\textbf{Proof of Proposition \ref{OUbougerol}} Let us consider a
second Ornstein-Uhlenbeck process $(U^{\lambda}_{t},t\geq 0)$
independent of the first one. Then, taking equation (\ref{OUeB})
for $U^{\lambda}_{t}$, we have:
\beq\label{eU}
    e^{\lambda t } U^{\lambda}_{t}= \delta_{( \frac{e^{2\lambda t}-1}{2 \lambda} )},
\eeq
where $(\delta_{t},t\geq 0)$ is a complex valued Brownian
motion starting from $z_{0}=1$. Thus:
\beq\label{TbU}
    T^{(\lambda)}_{b}(U^{\lambda})=\frac{1}{2\lambda}\ln \left(1+2\lambda T^{\delta}_{b}\right).
\eeq
Equation (\ref{thetaBMOU}) for $t=\frac{1}{2\lambda} \ln
\left(1+2\lambda T^{\delta}_{b}\right)$, equivalently:
$\alpha(t)=T^{\delta}_{b}$ becomes (we suppose that $z_{0}=1$):
\beqq
\theta^{Z^{\lambda}}_{T^{(\lambda)}_{b}(U^{\lambda}) } =
\theta^{Z^{\lambda}}_{\frac{1}{2\lambda}\ln \left(1+2\lambda
T^{\delta}_{b}\right) }= \theta^{\mathbb{B}}_{u=T^{\delta}_{b}}
\stackrel{(law)}{=} C_{a(b)}.
\eeqq
\begin{flushright}
  $\Box$
\end{flushright}
%%%%%%%%%%%%%%%%%%%%%%%%%%%%%%%%%%%%%%%%%%%%%%%%%%%%%%%%%%%%%%%%

%%%%%%%%%%%%%%%%%%%%%%%%%%%%%%%%%%%%%%%%%%%%%%%%%%%%%%%%%%%%%%%%
\subsection{On the distribution of $T^{\theta}_{-c,c}$ for an Ornstein-Uhlenbeck process}\label{OUgeneral}
%%%%%%%%%%%%%%%%%%%%%%%%%%%%%%%%%%%%%%%%%%%%%%%%%%%%%%%%%%%%%%%%
Now we turn to the study of the density function of:
\beqq
T^{(\lambda)}_{c} \equiv T^{\theta^{Z}}_{-c,c} \equiv \inf
\left\{t\geq 0 : \left|\theta^{Z}_{t}\right|=c \right\},
\eeqq
and its first moment.
%%%%%%%%%%%%%%%%%%%%%%%%%%%%%%%%%%%%%%%%%%%%%%%%%%%%%%%%%%%%%%%%
\begin{prop}\label{OUtime} Asymptotically for $\lambda$ large, for $z_{0}=1$, we
have:
\beq\label{ETchatasymp}
    2\lambda \: E\left[T^{(\lambda)}_{c}\right] - \ln\left( 2\lambda \right) \stackrel{\lambda\rightarrow \infty}{\longrightarrow} E\left[\ln\left(T^{\theta}_{-c,c}\right)\right],
\eeq
and:
\beq\label{ETchat2}
    E\left[\ln\left(T^{\theta}_{-c,c}\right)\right] = 2 \int^{\infty}_{0} \frac{dz}{\cosh \left(\frac{\pi z}{2}\right)} \ln\left(\sinh\left(cz\right)\right) +\ln\left(2\right) + c_{E},
\eeq
where $c_{E}$ is Euler's constant. \\
For $c<\frac{\pi}{8}$, we have the asymptotic equivalence:
\beq\label{ETchatasymp02}
    \frac{1}{\lambda}\left(E\left[T^{(\lambda)}_{c}\right] - E\left[\left(\sinh\left(B_{T^{\gamma}_{-c,c}}\right)\right)^{2}\right]\right) \stackrel{\lambda\rightarrow 0}{\longrightarrow} - \frac{1}{3}E\left[\left(\sinh\left(B_{T^{\gamma}_{-c,c}}\right)\right)^{4}\right].
\eeq
Equivalently:
\beq\label{ETchatderiv}
    \frac{d}{d\lambda} \Big|_{\lambda=0} E\left[T^{(\lambda)}_{c}\right] = \underset{\lambda\rightarrow 0}{\lim} \left[\frac{1}{\lambda}\left(E\left[T^{(\lambda)}_{c}\right] - E\left[T^{(0)}_{c}\right]\right)\right] = - \frac{1}{3}E\left[\left(\sinh\left(B_{T^{\gamma}_{-c,c}}\right)\right)^{4}\right]. \nonumber \\
\eeq
Moreover:
\beq\label{ETchatasymp0}
    E\left[\left(\sinh\left(B_{T^{\gamma}_{-c,c}}\right)\right)^{4}\right] = \int^{\infty}_{0} \frac{dz}{\cosh \left(\frac{\pi z}{2}\right)} \left(\sinh\left(cz\right)\right)^{4}.
\eeq
More precisely, for $c<\frac{\pi}{8}$:
\beq\label{ETchatasymp0cos}
    E\left[\left(\sinh\left(B_{T^{\gamma}_{-c,c}}\right)\right)^{4}\right] = \frac{1}{8} \left(\frac{1}{\cos(4c)}-4 \frac{1}{\cos(2c)} + 3\right),
\eeq
and asymptotically:
\beq\label{ETchatasymp0cosc0}
    E\left[\left(\sinh\left(B_{T^{\gamma}_{-c,c}}\right)\right)^{4}\right] \underset{c\rightarrow 0}{\simeq} 5c^{4}.
\eeq
\end{prop}
%%%%%%%%%%%%%%%%%%%%%%%%%%%%%%%%%%%%%%%%%%%%%%%%%%%%%%%%%%%%%%%%
%%%%%%%%%%%%%%%%%%%%%%%%%%%%%%%%%%%%%%%%%%%%%%%%%%%%%%%%%%%%%%%%
\textbf{Proof of Proposition \ref{OUtime}} \\
\textbf{\underline{$\lambda$ large}} \\
Let us return to equation (\ref{ETchat}). For
$\lambda\rightarrow\infty$, we have:
\beqq
    E\left[T^{(\lambda)}_{c}\right] &=& \frac{1}{2\lambda} E\left[\ln \left(1+2\lambda T^{\theta}_{-c,c}\right) \right] \\
                              &=& \frac{1}{2\lambda} E\left[\ln \left(2\lambda \left(T^{\theta}_{-c,c}+\frac{1}{2\lambda}\right) \right) \right] \\
                              &=& \frac{\ln\left( 2\lambda \right)}{2\lambda} + \frac{1}{2\lambda} E\left[\ln\left(T^{\theta}_{-c,c}+\frac{1}{2\lambda}\right)\right].
\eeqq
Thus:
\beqq 2\lambda \: E\left[T^{(\lambda)}_{c}\right] -
\ln\left( 2\lambda \right) \stackrel{\lambda\rightarrow
\infty}{\longrightarrow}
E\left[\ln\left(T^{\theta}_{-c,c}\right)\right],
\eeqq
which is precisely (\ref{ETchatasymp}). Moreover, by the integral
representation (\ref{logTcequation}) for
$E\left[\ln\left(T^{\theta}_{-c,c}\right)\right]$, we deduce
(\ref{ETchat2}). \\
\textbf{\underline{$\lambda$ small}} \\
We shall now study the case $\lambda\rightarrow 0$. We have that:
\beqq
    T^{(\lambda)}_{c} = \frac{1}{2\lambda}\ln \left(1+2\lambda T^{\theta}_{-c,c}\right).
\eeqq
For $c<\frac{\pi}{8}$, from Spitzer (\ref{spi}), (at least) the first two positive moments of
$T^{\theta}_{-c,c}$ are finite: $E\left[
\left(T^{\theta}_{-c,c}\right)^{p} \right]<\infty$,
$(p=1,2)$.
We make the elementary computation:
\beqq
    && \frac{1}{\lambda}\left(\frac{\ln \left(1+2\lambda x\right)}{2\lambda}-x\right)=\frac{1}{\lambda}\left(\frac{1}{2\lambda}\int^{1+2\lambda x}_{1}\frac{dy}{y}-x\right) \\
    && \stackrel{y=1+a}{=} \frac{1}{2\lambda^{2}}\int^{2\lambda x}_{0} \left(\frac{1}{1+a}-1\right) \; da
    \stackrel{a=2\lambda b}{=} -2\int^{x}_{0} \frac{b \; db}{1+2\lambda b} \stackrel{\lambda\rightarrow 0}{\longrightarrow} -x^{2}.
\eeqq
Consequently, by replacing $x=T^{\theta}_{-c,c}$, we have:
\beqq
    \frac{1}{\lambda} \left(E\left[ T^{(\lambda)}_{c} \right] - E\left[ T^{\theta}_{-c,c} \right] \right) =  E\left[ -2\int^{T^{\theta}_{-c,c}}_{0} \frac{b \; db}{1+2\lambda b} \right].
\eeqq
We may now use the dominated convergence theorem \cite{Bil78}, since the $(db)$ integral is majorized by $(T^{\theta}_{-c,c})^{2}$,
which is integrable. Thus:
\beqq
    \frac{1}{\lambda} \left(E\left[ T^{(\lambda)}_{c} \right] - E\left[ T^{\theta}_{-c,c} \right] \right) \stackrel{\lambda\rightarrow 0}{\longrightarrow} -E\left[ (T^{\theta}_{-c,c})^{2} \right].
\eeqq
Following the proof of Proposition \ref{logTc}, Bougerol's
identity and the scaling property yield:
\beqq
    E\left[\left(\sinh\left(B_{t}\right)\right)^{2}\right] &=& E\left[\left(\delta_{A_{t}}\right)^{2}\right] = E\left[A_{t} \left(\delta_{1}\right)^{2}\right] = E\left[A_{t}\right] \: E\left[\left(\delta_{1}\right)^{2}\right] \\
        &=& E\left[A_{t}\right] .
\eeqq
Thus, for $t=T^{\gamma}_{-c,c}$, we have:
\beqq
E\left[A_{T^{\gamma}_{-c,c}}\right] =
E\left[\left(\sinh\left(B_{T^{\gamma}_{-c,c}}\right)\right)^{2}\right].
\eeqq
Similarly:
\beqq
    E\left[\left(\sinh\left(B_{t}\right)\right)^{4}\right] &=& E\left[\left(\delta_{A_{t}}\right)^{4}\right] = E\left[\left(A_{t}\right)^{2}\left(\delta_{1}\right)^{4}\right] = E\left[\left(A_{t}\right)^{2}\right] \: E\left[\left(\delta_{1}\right)^{4}\right] \\
        &=& 3 E\left[\left(A_{t}\right)^{2}\right] .
\eeqq
Thus, for $t=T^{\gamma}_{-c,c}$, we have:
\beqq
    E\left[\left(A_{T^{\gamma}_{-c,c}}\right)^{2}\right] = \frac{1}{3} E\left[\left(\sinh\left(B_{T^{\gamma}_{-c,c}}\right)\right)^{4}\right].
\eeqq
So, because $A_{T^{\gamma}_{-c,c}}=T^{\theta}_{-c,c}$, we deduce (\ref{ETchatasymp02}). In order to prove
(\ref{ETchatderiv}), it suffices to remark that:
\beqq
E\left[T^{(0)}_{c} \right] = E\left[ T^{\theta}_{-c,c} \right] =
E\left[A_{T^{\gamma}_{-c,c}}\right] =
E\left[\left(\sinh\left(B_{T^{\gamma}_{-c,c}}\right)\right)^{2}\right] .
\eeqq
On the one hand, by using the density of
$B_{T^{\gamma}_{-c,c}}$:
\beqq
        E\left[\left(\sinh\left(B_{T^{\gamma}_{-c,c}}\right)\right)^{4}\right] &=& \int^{\infty}_{-\infty} \frac{dy}{2c} \: \frac{1}{\cosh(\frac{\pi y}{2c})} \left(\sinh\left(y\right)\right)^{4} \\
        &=& \int^{\infty}_{0} \frac{dy}{c} \: \frac{1}{\cosh(\frac{\pi y}{2c})} \left(\sinh y\right)^{4}  \\
        &\stackrel{z=\frac{y}{c}}{=}& \int^{\infty}_{0} dz \: \frac{1}{\cosh(\frac{\pi z}{2})} \left(\sinh (cz)\right)^{4},
\eeqq
which is finite if and only if $c<\frac{\pi}{8}$. In order
to prove this, it suffices to use the standard expressions:
$\sinh(x)=\frac{e^{x}-e^{-x}}{2}$ and
$\cosh(x)=\frac{e^{x}+e^{-x}}{2}$. On the other hand (note
$T\equiv T^{\gamma}_{-c,c}$), we remark that
$-B_{T}\stackrel{(law)}{=}B_{T}$ and \cite{ReY99}, ex.3.10,
$E\left[e^{kB_{T}}\right]=E\left[e^{\frac{k^{2}}{2}T}\right]=\frac{1}{\cos(kc)}$,
for $0\leq k<\pi(2c)^{-1}$, thus:
\beqq
    E\left[\left(\sinh\left(B_{T}\right)\right)^{4}\right] &=& \frac{1}{2^{4}}E\left[\left(e^{B_{T}}-e^{-B_{T}}\right)^{4}\right] \\
    &=& \frac{1}{2^{4}}E\left[ e^{4B_{T}} -4 e^{3B_{T}-B_{T}} + 6e^{2B_{T}-2B_{T}} -4 e^{B_{T}-3B_{T}} +e^{-4B_{T}} \right] \\
    &=& \frac{1}{2^{4}} \left(2E\left[ e^{4B_{T}}\right]- 8E\left[ e^{2B_{T}}\right]+6 \right) \\
    &=& \frac{1}{2^{3}} \left(\frac{1}{\cos(4c)}- 4\frac{1}{\cos(2c)}+3 \right),
\eeqq
which is precisely (\ref{ETchatasymp0cos}) and this is
finite if and only if $c<\frac{\pi}{8}$. Moreover, asymptotically
for $c\rightarrow 0$, by using the scaling property, we have:
\beqq
    E\left[\left(\sinh \left(B_{T^{\gamma}_{-c,c}}\right)\right)^{4}\right] &=& E\left[\left(\sinh \left(cB_{T^{\gamma}_{-1,1}}\right)\right)^{4}\right] \stackrel{c\rightarrow 0}{\simeq} c^{4} E\left[\left(B_{T^{\gamma}_{-1,1}}\right)^{4}\right] \\
    &=& c^{4} 3 E\left[\left(T^{\gamma}_{-1,1}\right)^{2}\right] = 5c^{4},
\eeqq
since $E\left[\left(T^{\gamma}_{-1,1}\right)^{2}\right]=5/3$
(see \cite{PiY03}; by using the notation of this
paper, Table 3: $E\left[X_{t}^{2}\right]=\frac{t(2+3t)}{3}$ for
$X_{t}=C_{1}$ and $t=1$). This asymptotics may also be obtained by
(\ref{ETchatasymp0cos}) by developing $\cos(4c)$ and $\cos(2c)$
into series up to the second order term and keeping the terms of
the order $c^{4}$.
\begin{flushright}
  $\Box$
\end{flushright}
%%%%%%%%%%%%%%%%%%%%%%%%%%%%%%%%%%%%%%%%%%%%%%%%%%%%%%%%%%%%%%%%
%%%%%%%%%%%%%%%%%%%%%%%%%%%%%%%%%%%%%%%%%%%%%%%%%%%%%%%%%%%%%%%%
\begin{rem}\label{OUz0D} If we slightly modify the above study for the
Ornstein-Uhlenbeck process by inserting a diffusion coefficient
$D$:
\beqq
Z_{t} = z_{0} + \sqrt{2D} \tilde{Z_{t}} - \lambda \int^{t}_{0} Z_{s} ds \: ,
\eeqq
we obtain:
\beq\label{OUz0DeB}
    Z_{t} &=& e^{-\lambda t} \left( z_{0} + \sqrt{2D} \int^{t}_{0} e^{\lambda s} d\tilde{Z_{s}} \right) \nonumber \\
                &=& e^{-\lambda t} \left( \mathbb{B}_{\alpha_{t}} \right),
\eeq
where in the second equation we used Dambis-Dubins-Schwarz
Theorem with
\beqq
    \alpha_{t}=2D \int^{t}_{0} e^{2\lambda s} ds = D \frac{e^{2\lambda t}-1}{\lambda}
\eeqq
\beqq \Rightarrow \alpha^{-1}_{t}= \frac{1}{2\lambda}  \ln
\left( 1+\frac{\lambda}{D} t\right).
\eeqq
Thus:
\beq\label{ETchatz0Dasymp}
    2\lambda \: E\left[T^{(\lambda)}_{c}\right] - \ln\left( \frac{\lambda}{D} \right) \stackrel{\lambda\rightarrow \infty}{\longrightarrow} E\left[\ln\left(T^{\theta}_{-c,c}\right)\right],
\eeq
because:
\beqq
    E\left[T^{(\lambda)}_{c}\right] &=& \frac{1}{2\lambda} E\left[\ln \left(1+\frac{\lambda}{D} T^{\theta}_{-c,c}\right) \right] \\
                              &=& \frac{1}{2\lambda} E\left[\ln \left(\frac{\lambda}{D} \left(T^{\theta}_{-c,c}+\frac{D}{\lambda}\right) \right) \right] \\
                              &=& \frac{\ln\left( \frac{\lambda}{D} \right)}{2\lambda} + \frac{1}{2\lambda} E\left[\ln\left(T^{\theta}_{-c,c}+\frac{D}{\lambda}\right)\right] .
\eeqq
Moreover:
\beq \label{ETchatz0D}
    E\left[\ln\left(T^{\theta}_{-c,c}\right)\right] &=& 2\ln(z_{0}) + E\left[\ln\left(T^{\theta^{(1)}}_{-c,c}\right)\right] \nonumber \\
        &=& 2\ln(z_{0}) + \int^{\infty}_{0} \frac{dz}{\cosh \left(\frac{\pi z}{2}\right)} \ln\left(\sinh\left(cz\right)\right) + \ln\left(2\right) + c_{E} , \nonumber \\
\eeq
where $T^{\theta^{(1)}}_{-c,c}$ denotes the first hitting time of the symmetric conic boundary of angle $c$ for a Brownian motion $Z$ starting from $1$. \\
For $\lambda$ small, we replace $2T^{\theta}_{-c,c}$ by
$\frac{z_{0}^{2}}{D}T^{\theta}_{-c,c}$ in the proof of Proposition
\ref{OUtime} ($\lambda$ small case) and we have:
\beqq
    T^{(\lambda)}_{c} &=& \frac{1}{2\lambda}\ln \left(1+\lambda \frac{z_{0}^{2}}{D} T^{\theta}_{-c,c}\right).
\eeqq
By repeating the previous calculation, we make the elementary computation:
\beqq
    && \frac{1}{\lambda}\left(\frac{\ln \left(1+\frac{z_{0}^{2}}{D} x\right)}{2\lambda}-\frac{z_{0}^{2}}{D}x\right)= -\frac{1}{2}\int^{x}_{0} \frac{\left(\frac{z_{0}^{2}}{D}\right)^{2}b \; db}{1+\lambda \frac{z_{0}^{2}}{D}b} \stackrel{\lambda\rightarrow 0}{\longrightarrow} -\left(\frac{z_{0}^{2}}{2D}\right)^{2}x^{2}.
\eeqq
We replace $x=T^{\theta}_{-c,c}$, and by the dominated convergence theorem \cite{Bil78}, for
$c<\frac{\pi}{8}$, we obtain:
\beqq
    \frac{1}{\lambda}\left(E\left[T^{(\lambda)}_{c}\right] - \frac{z_{0}^{2}}{2D} E\left[\left(\sinh\left(B_{T^{\gamma}_{-c,c}}\right)\right)^{2}\right]\right) &\stackrel{\lambda\rightarrow 0}{\longrightarrow}& - \frac{1}{3}\left(\frac{z_{0}^{2}}{2D}\right)^{2}E\left[(T^{\theta}_{-c,c})^{2}\right] \\
    &=& - \frac{1}{3}\left(\frac{z_{0}^{2}}{2D}\right)^{2}E\left[\left(\sinh\left(B_{T^{\gamma}_{-c,c}}\right)\right)^{4}\right],
\eeqq
where $E\left[\left(\sinh\left(B_{T^{\gamma}_{-c,c}}\right)\right)^{4}\right]$ is
given by (\ref{ETchatasymp0}), (\ref{ETchatasymp0cos}) and
asymptotically, for $c\rightarrow 0$ by (\ref{ETchatasymp0cosc0}).
\end{rem}
%%%%%%%%%%%%%%%%%%%%%%%%%%%%%%%%%%%%%%%%%%%%%%%%%%%%%%%%%%%%%%%%

%%%%%%%%%%%%%%%%%%%%%%%%%%%%%%%%%%%%%%%%%%%%%%%%%%%%%%%%%%%%%%%%
\noindent \textbf{Acknowledgements} \\
I would like to thank Alain Comtet for his contribution in the proof of Proposition \ref{density} and the Human Frontier Science Program for its support. This article represents the first part of my PhD thesis, under the supervision of Professors D. Holcman and M. Yor.
%%%%%%%%%%%%%%%%%%%%%%%%%%%%%%%%%%%%%%%%%%%%%%%%%%%%%%%%%%%%%%%%

%%%%%%%%%%%%%%%%%%%%%%%%%%%%%%%%%%%%%%%%%%%%%%%%%%%%%%%%%%%%%%%%

%%%%%%%%%%%%%%%%%%%%%%%%%%%%%%%%%%%%%%%%%%%%%%%%%%%%%%%%%%%%%%%%

\end{document}